\documentclass{article}
 \usepackage{latexsym}
 \usepackage{amsmath}
 \usepackage{amsfonts}
\usepackage{color}
\usepackage{amssymb}

\def\init{\setcounter{equation}{0}}

\newtheorem{theoreme}{Theorem }[section]
\newtheorem{proposition}[theoreme]{Proposition}
\newtheorem{lemma}[theoreme]{Lemma}
\newtheorem{definition}[theoreme]{Definition}

\newtheorem{remark}[theoreme]{Remark}

\newcommand{\beq}{\begin{equation}}
\newcommand{\eeq}{\end{equation}}
\newcommand\red{}

\def\bbbone{{\mathchoice {\rm 1\mskip-4mu l} {\rm 1\mskip-4mu l}
{\rm 1\mskip-4.5mu l} {\rm 1\mskip-5mu l}}}
\def\0{0}
\def\one{\bbbone}
\def\bel{\begin{lemma}}
\def\eel{\end{lemma}}
\def\bet{\begin{theoreme}}
\def\eet{\end{theoreme}}
\def\bed{\begin{definition}}
\def\eed{\end{definition}}
\def\bep{\begin{proposition}}
\def\eep{\end{proposition}}
\renewcommand\Re{\mathrm{Re}}

\def\rr{{\mathbb R}}
\def\zz{{\mathbb Z}}
\def\cc{{\mathbb C}}
\def\nn{{\mathbb N}}

\def\ii{\mathrm{ i}}

\def\s{{\rm s}}

\def\i{{\rm i}}

\def\e{{\rm e}}
\def\d{{\rm d}}

\def\C{{\cal C}}

\def\F{{\cal F}}

\def\S{{\cal S}}

\def\A{{\cal A}}

\def\Z{{\cal Z}}

\def\12{\frac{1}{2}}

\def\qed{\hfill$\Box$\medskip}
\def\proof{\noindent{\bf Proof.}\ \ }
\def\p{ \partial}

\def\cA{{\check A}}
\newcommand\cC{{\mathcal C}}
\def\cD{{\cal D}}

\def\cH{{\cal  H}}

\def\cF{{\cal F}}

\def\cS{{\cal S}}

\def\cA{{\cal A}}

\def\cZ{{\cal Z}}

\newcommand\hyp{\mathrm{hyp}}

\newcounter{smallarabics}
\newenvironment{arabicenumerate}
{\begin{list}{{\normalfont\textrm{(\arabic{smallarabics})}}}
  {\usecounter{smallarabics}\setlength{\itemindent}{0cm}
   \setlength{\leftmargin}{5ex}\setlength{\labelwidth}{4ex}
   \setlength{\topsep}{0.75\parsep}\setlength{\partopsep}{0ex}
   \setlength{\itemsep}{0ex}}}
{\end{list}}

\newcommand\ben{\begin{enumerate}}
  \newcommand\een{\end{enumerate}}

\begin{document}

\title{
A unified approach \\
to hypergeometric class functions}
 \author{J. Derezi\'{n}ski
\\
Department of Mathematical Methods in Physics,\\
Faculty of Physics, University of Warsaw, \\
Pasteura 5, 02-093 Warszawa, Poland\\
email jan.derezinski@fuw.edu.pl}
\maketitle

\begin{abstract} Hypergeometric class equations are given by second
  order differential operators 
in one  variable whose coefficient at the second derivative
  is a polynomial of degree $\leq2$, at the first derivative of degree
  $\leq1$ and the free term is a number.
Their solutions, called hypergeometric class functions,  include the Gauss hypergeometric function and its various
limiting cases. The paper presents a unified approach to
these functions. The main structure behind this approach is a
family of complex 4-dimensional Lie algebras, originally due to Willard
Miller. Hypergeometric class functions can be interpreted as 
eigenfunctions of the quadratic Casimir operator in a
representation of Miller's Lie algebra given by differential
operators in three complex variables. One obtains a unified
treatment of various properties of hypergeometric class functions
such as recurrence relations, discrete symmetries,  power series
expansions,  integral representations, generating functions and orthogonality of polynomial
solutions.
\end{abstract}

\tableofcontents

\init
\section{Introduction}
\label{s1}
This paper is devoted to the family of equations of the form
\beq\left(\sigma(z)\p_z^2+\tau(z)\p_z+
\eta\right) f(z)=0 ,
\label{req}\eeq
where $\sigma(z),\ \tau(z),\ \eta$ are polynomials with
\beq{\deg}\sigma\leq 2,\ 
{\deg}\tau\leq 1,\ \ \deg \eta=0.\label{req1}\eeq
Their solutions  include some of the most useful special functions
with many applications in physics and mathematics.

In the literature one can find several names 
for this family.
 In this paper, we will use the 
name {\em hypergeometric class equations} for equations of the form  \eqref{req} satisfying 
\eqref{req1}.\footnote{\red
In the book by Nikiforov-Uvarov
\cite{NU}, and also in \cite{De1,De2}, this family is called {\em equations of the
  hypergeometric type}. However,  Slavianov-Lay's book \cite{SL} suggests to
use the term ``type'' for smaller families, such as the families
(1)--(5) listed below. 
\cite{DIL} calls this family the  {\em grounded { Riemann} class}. The
name {\em {Riemann} class}, following the suggestions of
Slavianov-Lay,  is reserved in \cite{DIL}  for a somewhat wider
family, where $\eta(z)$ is allowed to be a rational function
such that $\eta\sigma$ is a polynomial of degree $\leq2$.
See Appendix \ref{Class} for more comments.}

Let us start with a short review of basic nontrivial
types of hypergeometric class equations. 
We will always  assume that $\sigma\neq0$.
Every class
 will be simplified
 by dividing by a constant and an affine change of the complex variable $z$.

\begin{arabicenumerate}

\item
\noindent
{\bf The ${}_2F_1$
or Gauss hypergeometric equation}

\[\left(z(1-z)\p_z^2+(c-(a+b+1)z)\p_z-ab\right)f(z)=0.\]

\item
\noindent
{\bf The ${}_2F_0$  equation}
\[\left(z^2\p_z^2+(-1+(1+a+b)z)\p_z+ab\right)f(z)=0.\]

\item
\noindent
{\bf The ${}_1F_1$ 
  or Kummer's  confluent equation}
 \[
(z\p_z^2+(c-z)\p_z-a)f(z)=0.\]

\item
\noindent
{\bf The ${}_0F_1$ equation}, closely related to the Bessel equation,
 \[
(z\p_z^2+c\p_z-1)f(z)=0.\]

\item
\noindent
{\bf The Hermite equation}
  \[
(\p_z^2-2 z\p_z-2a)f(z)=0.\]



\end{arabicenumerate}
 


 



In our work we collect various results about
 hypergeometric class
equations
 that can be stated and proven in a unified way, with as few restrictions on
 $\sigma$, $\tau$, $\eta$ as possible.
We believe that such an approach has a considerable
pedagogical and theoretical value. The pedagogical advantage of the unified
 approach is obvious: it reduces the need for repetitive arguments. From
 the theoretical point of view, in this way we  easily see the {\em
   coalescence} of various types. The properties  of  hypergeometric
 class equations that we describe depend {\em analytically} on the  the coefficients
$\sigma,\tau,\eta$. These properties include 
a pair of recurrence relations, a discrete symmetry, integral representations,
 power expansions around singular points, generating functions, the
 Rodriguez formulas for polynomials and their orthogonality.

The unified approach has its limitations. There are
some properties that are not easy to formulate in a unified way. For
instance, 
 only one pair of  recurrence relations 
depends analytically on the coefficients. If we restrict ourselves to  specific types,
then we often find a bigger number of
recurrence relations (eg. at least 12 in the case of the Gauss hypergeometric
equation). Another family of important properties not
included in our presentation are quadratic identities, which link
various types of hypergeometric class functions.
Thus this work should be compared to other studies of hypergeometric
class equations, such as \cite{De1},
where  specific types  are
described one by one.

The central structure behind the properties described in this paper is
a certain family of Lie algebras 
$m_{\alpha,\beta}$, described in Section \ref{Miller's Lie algebra}, which we propose to call  {\em Miller's
  Lie algebra}, since it was probably first introduced by Willard
Miller in \cite{M1}. Miller's Lie algebra
is defined as the span of
four elements    $N,A_+,A_-, \one$  satisfying the commutation relations
\beq\begin{array}{l}
[N,A_+]=A_+,\\[2mm]
[N,A_-]=-A_-,\\[2mm]
[A_+,A_-]=2\alpha N+\beta\one,
    \end{array}\label{nk1}\eeq
    where $\alpha,\beta$ are complex parameters.

Let $\sigma$ be a polynomial of degree $\leq2$ (as in \eqref{req}) and
$\kappa$ a polynomial of degree $\leq1$. Miller's Lie algebra with $\alpha=\frac{\sigma''}{2}$ and
$\beta=\kappa'$    can be represented by
    the following 1st order
differential operators acting on $\cc^3$:
\beq\begin{array}{l}
N:=t\p_t-s\p_s,\\[2mm]
A_+:=t\p_z+\sigma'(z)\p_s,\\[2mm]
A_-:=s\p_z+\sigma'(z)\p_t+\frac{\kappa(z)}{t}.
    \end{array}\label{n1.}\eeq
The operators \eqref{n1.}  can be 
restricted to functions
on the quadric \beq \sigma(z)-ts=0.\label{quad}\eeq

Miller's Lie algebra $m_{\alpha,\beta}$ commutes with
    \beq
\cC_{\alpha,\beta}:
=\12(A_-A_++A_+A_-)
+\alpha N^2+\beta N,
\label{casimir0}\eeq which
 will  be called the {\em  Casimir operator} of Miller's Lie algebra.
In the representation \eqref{n1.}, {\red the operator
  \eqref{casimir0}} restricted to the quadric
\ref{quad} and  to  the eigenspace $N=0$
is a 2nd order differential operator in $z$, which we denote
{\red $\cH(\sigma,\kappa)$.} Its  eigenvalue equation
\beq\big(\cH(\sigma,\kappa)+\omega\big)f=0\eeq
is precisely the hypergeometric class equation \eqref{req}, where \beq
\tau(z)=\kappa(z)+\sigma'(z),\quad\eta=\frac{\kappa'}{2}+\omega.\eeq

The operators $A_+$ and $A_-$ can be treated as ``root operators'' wrt
the ``Cartan algebra'' spanned by $N$. They lead to recurrence
relations for  ladders of solutions.
The involution of the quadric \eqref{quad} $(t,s,z)\mapsto(s,t,z)$ generates  a 
``Weyl symmetry'' of  Miller's Lie algebra. It leads to a discrete
symmetry of hypergeometric class equations.

We devote Section \ref{Basic properties of hypergeometric class
  equations} to the properties   of hypergeometric 
class equations and functions which follow directly from the
{\red action of root operators and Weyl symmetry} of Miller's Lie algebra: the {\em basic
  recurrence relations} and the {\em basic symmetry}.
Thanks to the recurrence relations, once we know a solution of  a
certain hypergeometric class equation, we know its solutions for a
whole ladder of parameters labelled by $n\in\zz$.

Certain ladders are special. One of them
contains
solutions which can be expressed in terms of elementary functions.
We call it the {\em Chebyshev
 ladder}, because it contains the well-known Chebyshev
polynomials. It is described in Subsection \ref{Chebyshev
  functions}. {\red See also \cite{DGR} Subsec. 3.3 and 4.5,  \cite{Poole},
Sec. 31.}

 Another ladder of solutions consists of polynomials. In the
 convention that we adopted, this ladder is descending---its elements
 are zero for $n>0$. We devote to this ladder the whole Section
 \ref{Hypergeometric  polynomials}.
The basic symmetry produces from the polynomial ladder  an ascending
ladder, which needs no separate discussion.

 For the Hermite equation the polynomial $\sigma$ does not have 
a zero.
For other types, that is for
equations reducible to the  ${}_2F_1$, ${}_1F_1$, ${}_2F_0$ and ${}_0F_1$ types, $\sigma$
has at least one zero.  We devote Section \ref{Singular point} to 
these equations, where without loss of generality we assume $\sigma(0)=0$. 
We can then write  a formal power series $F(\sigma,\kappa,\omega;z)$
which solves the  equation
\beq\big(\cH(\sigma,\kappa)+\omega\big)F=0,\eeq
normalized to be $1$ at
zero.
This power series is convergent if $\sigma'(0)\neq0$, that is for
${}_2F_1$, ${}_1F_1$, and ${}_0F_1$. If
$\sigma'(0)=0$, that is for ${}_2F_0$, {\red it either terminates and
  is a polynomial, or  is divergent. Nevertheless, even in the
  divergent case this series  is
asymptotic to  a well-defined
function ${}_2F_0$, which we separately discuss in Appendix \ref{appendix}.} Thus we obtain a unification of 4 types of hypergeometric
functions in a single function $F(\sigma,\kappa,\omega;z)$, which we  call the {\em unified
  hypergeometric function}. It depends meromorphically on $5$ complex
parameters (two parameters for $\sigma$, two for $\kappa$ and 
one for $\omega$).

If $\sigma(0)=0$ and $\sigma'(0)\neq0$, it is natural to assume
$\sigma'(0)=1$ and to normalize the
unified hypergeometric function
by dividing it by $\Gamma(m+1)$, where
$m=\frac{\kappa(0)}{\sigma'(0)}$.
This is sometimes
called {\em Olver's normalization}. Thus we obtain the {\em Olver
  normalized unified hypergeometric function}
$\mathbf{F}(\sigma,\kappa,\omega;z)$, which depends analytically on $4$
complex parameters (one parameter for $\sigma$, two for $\kappa$ and 
one for $\omega$).

If $\sigma(0)=0$ and $\sigma'(0)\neq0$, the hypergeometric class  has an additional discrete
symmetry, {\red which we call  the {\em power symmetry}. This transformation
involves gauging with the function $x^m$, that is,
  $(\cdot)\mapsto x^m(\cdot) x^{-m}$. (In the literature this
transformation is sometimes called an {\em F-homotopy}.)}

If $\sigma(0)=0$, and $\sigma''\neq0$ or $\kappa'\neq0$,
that is for the  ${}_2F_1$, ${}_1F_1$ and  ${}_2F_0$  equations,
we have yet another discrete symmetry, which we call the {\em
  inversion symmetry}. For the ${}_2F_1$ equation it involves the
change of the variable  $z\mapsto \frac1z$ and the
transformed equation is still of the form
${}_2F_1$. In the confluent cases
it interchanges the ${}_1F_1$ and  ${}_2F_0$  equations, and
it involves $z\to-\frac1z$. It is however not defined for the ${}_0F_1$
equation.

In Section \ref{Integral representations} we describe a unified approach
to 
representations
of hypergeometric class functions in terms of integrals of elementary functions. 
The ${}_2F_1$,  ${}_1F_1$, ${}_2F_0$ and Hermite functions
can be represented as {\em Euler integrals}, that
is,
\beq \int_\gamma p(s)(s-z)^{-n-1}\d s,\eeq
where $p$ is an elementary function, $n\in\cc$  and $\gamma$ is an
appropriate contour.
${}_1F_1$, ${}_0F_1$ and Hermite functions
possess   representations in the form of the {\em Laplace integral}
\beq \int_\gamma q(s)\e^{zs}\d s,\eeq
where $q$ is an elementary function and $\gamma$ is an
appropriate
contour.

As we mentioned above, some ladders of solutions of hypergeometric
class equations consist of polynomials. These polynomials include the
famous {\em classical orthogonal polynomials} (sometimes called the
{\em very classical orthogonal polynomials}).
They allow for a very elegant unified treatment, which includes not
only recurrence relations, but also 
generating functions, the famous Rodrigues formula and the
orthogonality relations. All of this is briefly described in Section \ref{Hypergeometric  polynomials}.

As we mentioned above,  if  $\sigma(0)=0$ and $\sigma'(0)\neq0$, that is in the cases
equivalent to
 types
${}_2F_1$, ${}_1F_1$, and
${}_0F_1$, we have the power symmetry. Because of that,
these equations have two linearly
independent
solutions with a distinct behavior at
zero: the unified hypergeometric function, which is analytic at zero,
and another solution behaving as $\sim z^m$, where
$m=\frac{\kappa(0)}{\sigma'(0)}$. The linear independence  breaks
down for $m\in\zz$, when both solutions are proportional to one
another. This case is called {\em degenerate} and is discussed in Section
\ref{Degenerate case}. The Olver normalizzed unified
hypergeometric function in this case has an additional
integral representation and an elegant generating function.

The present work should be used 
only as an ``invitation'' to  hypergeometric class 
functions. It leaves out many of their properties, which are difficult to describe in a
unified way.
For instance, as we mentioned above, various types of hypergeometric class functions
possess additional recurrence relations and additional
symmetries. However only those derived directly from  Miller's Lie algebra can
be described in a fully unified fashion.

Almost all our discussion is algebraic, without functional
analysis. We are aware that it is natural to view
 hypergeometric class operators as closed (or self-adjoint) Sturm-Liouville operators
on appropriate
weighted $L^2$ spaces.
It would be also interesting to
consider representations of Miller's Lie algebra in Hilbert
spaces. This would, however, require breaking
our discussion into separate types.
The only place where we use some elements of functional analysis
is Section  \ref{Hypergeometric  polynomials}
 about
classical orthogonal polynomials. We show how to view them as
eigenfunctions of certain self-adjoint Sturm-Liouville operators on 
weighted $L^2$ spaces with appropriate boundary conditions---this can
be done in a rather unified fashion.

The literature on hypergeometric class equations is
very large. Usually, each type is considered separately, without an
attempt of a unified treatment. Let us list some of the more famous
treatises about these equations: \cite{EMOT,Ho,MOS,MF,R,Wa,WW}

Under the name ``equations of the
hypergeometric type'' they were considered in a unified way in the book by 
Nikiforov--Uvarov \cite{NU}. This book was one of the two   main inspirations for
our article. A part of the  material of our work is
adapted from \cite{NU}, notably the material of Subsections 
\ref{Ladders of solutions},  \ref{s.9}, \ref{s.12},
\ref{s.13}.

We use also ideas of Miller, notably
in Section \ref{Miller's Lie algebra}. As we mentioned above, Miller's
Lie algebra was defined in \cite{M1}. Miller was an early champion of
the use of Lie algebras in the theory of special functions. 
To my 
knowledge, he was the first to note that various types of 
hypergeometric class admit a larger symmetry algebra. This topic was 
further developed in \cite{De2}.

{\red There are many works that contain elements of a unified theory of
hypergeometric class equations,  the idea of  recurrence
relations and factorizations. Among early ones let us mention the
classic work of Infeld and Hull \cite{IH}
and of Truesdell \cite{Tr}. Later treatments include \cite{CKS,SHD}.}

More complete analysis of various types hypergeometric class
equations, including the Lie-algebraic interpretation of their
recurrence relations and discrete symmetries, can be found in the
literature, notably in the works of Miller, and also in
\cite{De1,De2}.

The family of hypergeometric class polynomials that form an orthogonal basis in an
appropriate weighted Hilbert space consists essentially of Jacobi, Laguerre and
Hermite polynomials, often called
 {\em classical 
orthogonal polynomials}.
It has an especially large literature, e.g.
\cite{NU,R}. In the more recent literature
the name ``classical orthogonal polynomials'' is sometimes given to a
broader family, given by the so-called Askey scheme. Some authors
proposed to use the name {\em very classical orthogonal polynomials}
for
the family consisting of Jacobi, Laguerre and
Hermite polynomials.

\init
\section{Miller's Lie algebra}
\label{Miller's Lie algebra}

In  this section we introduce a certain two-parameter family of 4-dimensional Lie
algebras, which to our knowledge was first introduced by Willard Miller. We
will call it {\em Miller's Lie algebra} and denote by $m_{\alpha,\beta}$
with $\alpha,\beta\in\cc$. We will also describe its
{\em Casimir operator} $\cC_{\alpha,\beta}$, that is, a quadratic
expression in elements of
$m_{\alpha,\beta}$, which in any representation
commutes with the whole Lie algebra. In other words,
$\cC_{\alpha,\beta}$ belongs to the center of the enveloping
algebra of $m_{\alpha,\beta}$. We will also describe a family of representations
of $m_{\alpha,\beta}$
by 1st order differential operators on certain 2nd degree surfaces in
$\cc^3$.

In appropriate coordinates the eigenvalue equation for the Casimir operator 
equation will have the form of a
hypergeometric class equation. A number of properties
of hypergeometric class
equations will have a simple  interpretation in terms of
properties of Miller's Lie algebra. They include
the basic symmetry, basic 
factorizations and basic recurrence relations. They will be described in the next section.

\subsection{Three low-dimensional Lie algebras}

Let us first introduce three low-dimensional complex  Lie algebras. 

\begin{enumerate}
  \item $sl(2,\cc)$. It is spanned by $N,A_+,A_-$ satisfying the 
    commutation relations 
\beq\begin{array}{l}
[N,A_\pm]=\pm A_\pm,\\[2mm]
[A_+,A_-]=2 N. 
    \end{array}\label{sl2c}\eeq
The  operator
    \begin{align}
\cC:=&A_-A_++ N(N+1)
\label{casimir.1.}\\
=&A_+A_-+ N(N-1)
\label{casimir.2.}\\
=&\12(A_-A_++A_+A_-)
+ N^2
\label{casimir.}\end{align}
commutes with all elements  of
$sl(2,\cc)$. It is often called the {\em Casimir (operator)}.

Here is a typical representation of $sl(2,\cc)$:
 \begin{align}
   N&:=x\partial_x-y\partial_y,\\
   A_+&:=x\partial_y,\quad    A_-:=y\partial_x. 
        \end{align}
  \item The so-called {\em complex oscillator Lie algebra}, denoted $osc(\cc)$. 
    It is spanned by $N,A_+,A_-,\one$ satisfying the 
    commutation relations 
\beq\begin{array}{l}
[N,A_\pm]=\pm A_\pm,\\[2mm]
      [A_+,A_-]=\one,\\[2mm]
      [\one,A_\pm]=[\one,N]=0. 
    \end{array}\label{osc}\eeq
The following operator commutes with $osc(\cc)$, and will be 
called the {\em Casimir} of $osc(\cc)$:
    \begin{align}
\cC:=&A_-A_+
+N+\tfrac12\label{casimir.1a}\\
=&A_+A_-
+N-\tfrac12\label{casimir.2a}\\
=&\12(A_-A_++A_+A_-)+ N.
   \label{casimira}\end{align}

 Here is a typical representation of $osc(\cc)$: 
 \begin{align}
   N&:=\frac12(-\partial_x^2+x^2),\\
   A_\pm&:=\frac1{\sqrt2}(x\mp\partial_x).
        \end{align}
 Note that $N$ is the quantum harmonic oscillator, which justifies the
 name of this Lie algebra.

  \item The Lie algebra of {\em Euclidean movements of the plane}, denoted
    $\cc^2\rtimes so(2,\cc)$. 
    It is spanned by $N,A_+,A_-$ satisfying the 
    commutation relations 
\beq\begin{array}{l}
[N,A_\pm]=\pm A_\pm,\\[2mm]
      [A_+,A_-]=0.
\end{array}\label{euc}\eeq
\end{enumerate}
The following operator commutes with $\cc^2\rtimes so(2,\cc)$, and will be 
called the {\em Casimir} of $\cc^2\rtimes so(2,\cc)$:
\begin{align}
\cC:=&A_-A_+
\label{casimirb}\end{align}

Here is a typical representation of     $\cc^2\rtimes so(2,\cc)$:
 \begin{align}
   N&:=x\partial_y-y\partial_x,\\
   A_\pm&:=\partial_x\pm\i\partial_y.
        \end{align}
        Thus $N$ is the generator of rotations of the plane and
        $A_\pm$ generate translations.

\subsection{The family of Lie algebras introduced by W. Miller}

The three Lie algebras introduced in the previous subsection
can be joined in a single family depending on
two complex parameters $\alpha,\beta$. This family will be denoted by
$m_{\alpha,\beta}$. It was introduced by Willard Miller in  
\cite{M1}. We will call it  {\em Miller's Lie algebra}.

$m_{\alpha,\beta}$ is defined as the complex Lie algebra spanned 
by $N,A_+,A_-,\one$ satisfying the 
    commutation relations 
\beq\begin{array}{l}
[N,A_\pm]=\pm A_\pm,\\[2mm]
      [A_+,A_-]=2\alpha N+\beta\one,\\[2mm]
      [\one,A_\pm]=[\one,N]=0. 
\end{array}\label{eucl}\eeq
It is easy to see that
\begin{align}\label{miller1}
  m_{\alpha,\beta}&\simeq sl(2,\cc)\oplus\cc,&&\quad \alpha\neq0;\\\label{miller2}
                     m_{\alpha,\beta}&\simeq osc(\cc),&&\quad
                                      \alpha=0,\quad\beta\neq0;\\\label{miller3}
                                                           m_{\alpha,\beta}&\simeq
                                                                            \cc^2{\rtimes}
                                                                            so(2,\cc)\oplus\cc,&&\quad
                                                                            \alpha=0,\quad\beta=0.
\end{align}

Define the linear map ${\red \varepsilon}: m_{\alpha,\beta}\to m_{\alpha,-\beta}$ by 
\begin{align}
  {\red \varepsilon}(N):=-N,\quad  
  {\red \varepsilon}(A_\pm):=A_\mp,\quad  
  {\red \varepsilon}(\one):=\one. 
\label{tau}\end{align} 
Then ${\red \varepsilon}$ is an isomorphism. The identity automorphism together with
${\red \varepsilon}$ form a group isomorphic to $\zz_2$, which
will be called the {\em Weyl group} of $m_{\alpha,\beta}$.

Similarly, for $n\in\cc$, the linear map $\pi_n:m_{\alpha,\beta}\to
m_{\alpha,\beta-2n\alpha}$ given by
\begin{align}
  \pi_n(N):=N+n\one,\quad 
  \pi_n(A_\pm):=A_\pm,\quad 
  \pi_n(\one):=\one, 
\label{pin}\end{align}
is an isomorphism.

\subsection{Casimir}

Consider a representation of
$m_{\alpha,\beta}$ on a vector space $\cZ$.
Following Miller \cite{M1},
introduce the following operator $\cC$, called the {\em Casimir of $m_{\alpha,\beta}$}:
\begin{align}
\cC_{\alpha,\beta}={}\cC:=&A_-A_++\alpha N(N+1)
+\beta(N+\tfrac12)\label{casimir.1}\\
=&A_+A_-+\alpha N(N-1)
+\beta(N-\tfrac12)\label{casimir.2}\\
=&\12(A_-A_++A_+A_-)
+\alpha N^2+\beta N.
\label{casimir}\end{align}
As Miller noted, $\cC$ commutes with the whole 
Lie algebra: 
\begin{align}
  N\cC&=\cC N,\qquad
  A_\pm \cC=\cC A_\pm.\label{commi}
\end{align}

Extend the isomorphisms ${\red \varepsilon}$ and $\pi_n$ defined in \eqref{tau} and
\eqref{pin} to the algebra of operators on $\cZ$.
Then
\begin{align}{\red \varepsilon}(\cC_{\alpha,\beta})&=\cC_{\alpha,-\beta},\\
\pi_n(\cC_{\alpha,\beta})&=\cC_{\alpha,\beta+2n\alpha}+(\alpha n^2+\beta n)\one.
\end{align}

\subsection{Ladders}

By a {\em two-sided ladder} we mean a subset of $\cc$ of the form
 $n_0+\zz$, where
 $n_0\in\cc$.
A subset of the form $n_0+\nn_0$ or $n_0-\nn_0$ will be called a
{\em one-sided ladder}  ({\em ascending} or {\em descending}). A
subset of the form $\{n_0,n_0+1,\dots,n_0+n\}$ for some $n\in\nn_0$ will
be called a {\em finite ladder}.

It is easy to see that 
 if  {\red the representation of $m_{\alpha,\beta}$ on }$\cZ$ is irreducible and $N$ possesses an eigenvalue, then the spectrum of $N$ is
 a ladder. In fact this follows from
\beq Nv=nv \quad\Rightarrow\quad NA_\pm v=(n\pm1)A_\pm v.\eeq

For $n\in\cc$ we define
\beq \cZ^n:=\{v\in\cZ\ :\quad\quad Nv=nv\}.\label{zetn}\eeq
Let $\cC^n$, resp. $A_\pm^{n\pm\frac12}$, be the operator $\cC$, 
resp. $A_\pm$, restricted to $\cZ^n$. \eqref{commi} can be
rewritten as
\beq
  A_\pm^{n\pm\frac12} \cC^n=\cC^{n\red\pm1} A_\pm^{n\pm\frac12}.\label{commi1}
\eeq

An irreducible representation with a one-sided or finite ladder with the lowest,
resp. highest
element $0$ will
be called a {\em lowest}, resp. {\em heighest weight
  representation}. It follows from \eqref{casimir.1},
resp. \eqref{casimir.2} that
\begin{align}\cC^{\red0}=\frac{\beta}2&\quad \text{for highest weight
                               representations},\\
  \cC^{\red0}=-\frac{\beta}2&\quad \text{for lowest weight representations}.
                    \end{align}

\subsection{Representation by 1st order differential operators}
\label{s.5}

Let $\sigma$ be a polynomial of degree $\leq2$ and $\kappa$ of degree $\leq1$.
Consider $\cc^3$ with elements denoted by $(t,s,z)$. 
Define the operators 
\beq\begin{array}{l}
N:=t\p_t-s\p_s,\\[2mm]
A_+:=t\p_z+\sigma'(z)\p_s,\\[2mm]
A_-:=s\p_z+\sigma'(z)\p_t+\frac{\kappa(z)}{t}.
    \end{array}\label{n1}\eeq
    $A_+, A_-, N,\one$ span a Lie algebra
with commutation relations
\beq\begin{array}{l}
[N,A_+]=A_+,\\[2mm]
[N,A_-]=-A_-,\\[2mm]
[A_+,A_-]=\sigma'' N+\kappa'\one.  
\end{array}\label{nk1}\eeq  
{\red Thus
it }
 is  a representation of  Miller's  Lie algebra 
$m_{\alpha,\beta}$ with \beq\alpha=\frac{\sigma''}{2},\quad\beta=\kappa'.\eeq

The operators {\red\eqref{n1}} commute with the multiplication by
$\sigma(z)-ts$. Therefore, we can restrict them to
analytic functions on $\Omega$, the universal cover of the manifold
\beq
\big\{(s,t,z)\in\cc^3\ :\ \sigma(z)-ts=0,\ t\neq0,\ s\neq 0\big\}.
\label{n2}\eeq
Let  $\A(\Omega)$ denote the space of analytic functions on
 $\Omega$.
 The Lie algebra $m_{\alpha,\beta}$ represented by \eqref{nk1}
acting on $\cA(\Omega)$ 
will be denoted $m(\sigma,\kappa)$.

By \eqref{casimir}, the Casimir of $m(\sigma,\kappa)$ is following
differential operator
on $\cA(\Omega)$:
\beq\begin{array}{rl}
\cC=&st\p_z^2+\sigma'(z)(s\p_s+t\p_t+1)\p_z 
+\kappa(z)\p_z\\[2mm]
&+\frac{\sigma''}{2}
\left((t\p_t)^2+(s\p_s)^2-2t\p_ts\p_s+t\p_t+s\p_s\right)\\[3mm]&
+\frac{(\sigma'(z))^2}{ts}t\p_ts\p_s+
\frac{\kappa(z)\sigma'(z)}{ts}s\p_s+\kappa'(t\p_t-s\p_s+\12). 
\end{array}\label{faf}\eeq 
Clearly, $\cC$  commutes with $m(\sigma,\kappa)$.

Introduce new coordinates on $\cc^3$:
\beq v=\sigma(z)-ts,\ \ \       \ w=t,\ \ \  \ \check z=z.\label{coor}
\eeq
The inverse transformation is
\[\begin{array}{l}
t=w,\ \ \ \ \ s=\frac{\sigma(z)-v}{w},\ \ \ \ \ z=\check z.\end{array}\]
We have
\[\begin{array}{l}
\p_z=\p_{\check z}+\sigma'(z)\p_v,\ \ \ \ \ 
\p_t=\p_{w}-s\p_v,\ \ \ \ \ 
\p_s=-t\p_v
\end{array}\]

Clearly, \eqref{n2} is $v=0$. 
Therefore the operators (\ref{n1}) in new coordinates, after restricting to
the surface
(\ref{n2}) (and dropping ``checks''), are
\beq\begin{array}{l}
{} N=w\p_w,\\[2mm]
{} A_+=w\p_z,\\[2mm]
{} A_-= 
\frac{1}{w}\left(\sigma(z)\p_z+\sigma'(z)w\p_w+\kappa(z)\right).
\end{array}\label{n1a}\eeq
The Casimir operator is
\[\begin{array}{ll}
{}\C&={} A_-{} A_++\frac{\sigma''}{2}{} N({} N+1)
+\kappa'{} (N+\12)\\[2mm]
&=\sigma(z)\p_z^2+(\kappa(z)+\sigma'(z)(w\p_w+1))\p_z
+\frac{\sigma''}{2}w\p_w(w\p_w+1)+\kappa'(w\p_w+\12).
\end{array}\]

Consistently with \eqref{zetn}, for $n\in\cc$ let
\[\Z^n:=\{g\in\A(\Omega)\ :\ \quad Ng=ng\}.\]
As in the previous subsection, $\cC^n$, resp. $A_\pm^{n\pm\frac12}$ denote the operator $\cC$, 
resp. $A_\pm$ restricted to $\cZ^n$. 
Clearly, in the coordinates (\ref{coor}) 
\[\Z^n=\{w^nF(z)\ :\quad F\in\A(\Theta)\},\]
where $\Theta$ is the universal covering of $\cc\backslash\{\hbox{zeros
  of }\sigma\}$.

 We have
 \beq\label{casimir1}
 \C^n+\omega:=
\sigma(z)\p_z^2+(\kappa(z)+\sigma'(z)(n+1))\p_z
                                +\frac{\sigma''}{2}n(n+1)+\kappa'(n+\tfrac12)+\omega.\eeq
\[A_{+}^{n+\frac12}:=\p_z,\ \ \ \ A_{-}^{n-\frac12}:=\sigma(z)\p_z+\kappa(z)+n\sigma'(z).\]


\subsection{Implementation of isomorphisms}
Let us go back to $\cc^3$ {\red in the coordinates $s,t,z$.}
{\red An} implementation of the isomorphism $\pi_n$ defined in  
\eqref{pin} is simple:
\beq
\pi_ng(\cdot)=t^{-n}\cdot t^n.\eeq

Let us now describe an implementation of the isomorphism ${\red \varepsilon}$.
Let $\rho(z)$ solve
\beq(\sigma(z)\p_z-\kappa(z))\rho(z)=0,\label{e1b.}\eeq
(which defines $\rho(z)$ up to a coefficient).
We introduce the following transformation on $\A(\Omega)$:
\[Tg(t,s,z):=\rho(z)g(s,t,z).\]

\bet\label{implem}
\beq TNT^{-1}=-N,\label{sys1}\eeq
\beq TA_+T^{-1}=A_-,\label{sys2}\eeq
\beq TA_-T^{-1}=A_+.\label{sys3}\eeq
\beq T\C T^{-1}=\C.\label{sys4}\eeq
where on the left we have the operators from $m(\sigma,
\kappa)$ and on the right from $m(\sigma,-\kappa)$.
\eet

\proof \eqref{sys1} is immediate.

Consider $g(t,s,z)\in\A(\Omega)$.
Then
\[\begin{array}{rl}
T^{-1} g(t,s,z)&=\rho^{-1}(z)g(t,s,z),\\[3mm]
A_+T^{-1}
g(t,s,z)&=(t\rho^{-1}(z)\p_z-\frac{\kappa(z)}{\sigma(z)}
t\rho^{-1}(z)+\sigma'(z) \rho^{-1}(z)\p_s)g(s,t,z)\\[3mm]
&=(t\rho^{-1}(z)\p_z-\frac{\kappa(z)}{s}
\rho^{-1}(z)+\sigma'(z) \rho^{-1}(z)\p_s)g(s,t,z)
\\[3mm]
T A_+T^{-1}
g(t,s,z)
&=(s\p_z-\frac{\kappa(z)}{t}
+\sigma'(z)\p_t)g(t,s,z)
.
\end{array}\] 
This proves (\ref{sys2}).
The proof of (\ref{sys3}) is similar.

To show (\ref{sys4}) we use 
 the second line of (\ref{faf}), (\ref{sys1}), (\ref{sys2}) and (\ref{sys3}).
\qed

In the coordinates \eqref{coor}
the  symmetry equals
\[\begin{array}{l}
    T
    g(w,z)=\rho(z)
    g(\frac{\sigma(z)}{w},z)
\end{array}.\]
Let $T^n$ be $T$ restricted to $\cZ^n$.
Clearly,
\beq \cC^n T^n=\cC^{-n} T^{-n}.\eeq
We easily see that in the coordinates \eqref{coor}, {\red for all $F\in\cA(\Omega)$,}
\[T^n F(z)=\rho(z)\sigma^n(z)F(z).\]


\init
\section{Basic properties of hypergeometric class equations}
\label{Basic properties of hypergeometric class equations}

In this section we introduce hypergeometric class
equations. They will be presented as the eigenvalue equation of a
certain second order differential operator $\cH(\sigma,\kappa)$. We
discuss a number of  properties of these equations
that can be
described in a unified way:
the basic symmetry, basic 
factorization and basic recurrence relations.

{\red The operator $\cH(\sigma,\kappa)$ is essentially the Casimir $\cC$} of Miller's
Lie algebra at $N=0$, and all its properties discussed in this section
follow from the properties of
Miller's Lie algebra analyzed  in the previous section.
However,  this section can be read independently.

\subsection{Remark about notation}

Let  $a,b,c$ be  complex functions. 
Instead of saying that we consider the equation 
\beq 
\Big(a(z)\p_z^2+b(z)\p_z+
c(z) \Big)f(z)=0,\eeq 
we will usually say that the equation is {\em given by the operator} 
\beq\cA:=a(z)\p_z^2+b(z)\p_z+
c(z) 
\label{reqs}\eeq

Instead of $\cA$, we will sometimes use the 
notation $\cA(z,\p_z)$ to indicate the variable that is used in the
given operator. This is useful when we consider a change of variables.

\subsection{Parametrization of hypergeometric class operators}
\label{s.2}

As was described in the introduction, the main topic of the paper are
equations
given by operators of the form
\beq\sigma(z)\p_z^2+\tau(z)\p_z+
\eta,
\label{req0}\eeq
where $\sigma,\tau,\eta$ satisfy
the conditions \eqref{req1}. However, the parametrization of these
equations with $\sigma,\tau,\eta$ is not always convenient. {\red We will
usually prefer to use a different parametrization,  described below.}

Let $\sigma(z),\ \kappa(z),\omega$ be polynomials with
\beq{\deg}\sigma\leq 2,\ 
{\deg}\kappa\leq 1,\ \ \deg \omega=0\quad\text{(in other words $\omega\in\cc$)}.\label{war}\eeq
Let us define the  differential operator
\begin{align}     \label{coro}
{\red \cH}(\sigma,\kappa)&:=
\sigma(z)\p_z^2+\big(\sigma'(z)+\kappa(z)\big)\p_z+\tfrac12\kappa'
\\
&=\p_z\sigma(z)\p_z+\tfrac12\big(\p_z\kappa(z)+\kappa(z)\p_z\big).
\notag
\end{align}
Clearly, the class of operators (\ref{req0}) coincides with the class of operators
\beq
{\red \cH}(\sigma,\kappa)+\omega
\label{e1}\eeq

Note that ${\red \cH}(\sigma,\kappa)$ coincides with ${\red \cC}^0$,
the Casimir operator for Miller's Lie
algebra restricted to the subspace $N=0$, {\red expressed in the variable
$z$, } see
\eqref{casimir1}.



\subsection{Basic symmetry}
\label{s.3}

Recall from \eqref{e1b.} that $\rho(z)$ is defined as a solution of
\beq(\sigma(z)\p_z-\kappa(z))\rho(z)=0.\label{e1b}\eeq
We have the identity
\beq
{\red \cH}(\sigma,\kappa)=
\rho^{-1}(z)\p_z\sigma(z)\rho(z)\p_z+\tfrac12\kappa'.
\label{e1a}\eeq

The following theorem describes a certain symmetry of the entire family
of hypergeometric class
equations. We will call it the {\em basic symmetry}.

\bet
We have
\beq\begin{array}{l}
\rho(z){\red \cH}(\sigma,\kappa)\rho^{-1}(z)
={\red \cH}(\sigma,-\kappa).
\end{array}
\label{asa}\eeq
Hence,
\beq
\big({\red \cH}(\sigma,\kappa){\red +\omega\big)}F=0\ \Rightarrow\
  \big({\red \cH}(\sigma,-\kappa){\red+\omega\big)}\rho F=0.
\label{syme0}\eeq\label{syme}
\eet

\proof Using $\sigma(z)\partial_z\rho^{-1}(z)=-\kappa(z)\rho^{-1}(z)$,
we obtain
{\red\begin{align}
  \rho(z)\cH(\sigma,\kappa)\rho^{-1}(z)&=\partial_z\sigma(z)\rho(z)\partial_z\rho^{-1}(z)+\frac12\kappa'\\
=\sigma(z)\partial_z^2+(\sigma'(z)-\kappa(z))\partial_z-\frac12\kappa'&=\cH(\sigma,-\kappa).\end{align}}
\qed

 Applying twice the basic  symmetry   we
get the identity.
Hence we obtain a group of symmetries of {\red the hypergeometric class}
isomorphic to $\zz_2$.

Note that the basic symmetry of Theorem \ref{syme} corresponds to the
symmetry ${\red \varepsilon}$ of Miller's Lie algebra, see \eqref{tau} and Theorem \ref{implem}.



\subsection{Basic pair of factorizations}
\label{s.4}

It is often useful to represent a 2nd order operator as
a product of two 1st order operators plus a constant.
In this section we describe a pair of  such factorizations of
hypergeometric class
operators, which can
be easily formulated in a unified way. 
It is convenient to formulate them for a family  indexed by $n\in\cc$.
These factorizations lead 
to recurrence relations for hypergeometric class functions.

Fix $\sigma$, a polynomial of degree $\leq2$, $\kappa_0$, a polynomial
of degree $\leq1$ and $\omega_0\in\cc$.
  We set
\begin{align}
  \kappa_n(z):=&n\sigma'(z)+\kappa_0(z),\label{kappan}\\
  \omega_n:=&n^2\frac{\sigma''}{2}+n\kappa_0'+\omega_0.\label{lambdan}
               \end{align}
  {\red Note that $\cH(\sigma,\kappa_n)+\omega_n$ coincides with 
  $\cC^n+\omega_0$, where $\cC^n$ is the Casimir for $\sigma,\kappa_0$
  restricted to $N=n$, see \eqref{casimir1}.}

\bet
(1) Factorization properties
\[\begin{array}{rl}
{\red \cH}\big(\sigma,\kappa_n\big)+\omega_n&=
\Big(\sigma(z)\p_z+\kappa_{n+1}(z)\Big)\p_z+n(n+1)\frac{\sigma''}{2}
+(n+\12)\kappa_0'+\omega_0\\[3mm]
&=\p_z\Big(\sigma(z)\p_z+\kappa_n(z)\Big)
+n(n-1)\frac{\sigma''}{2}+(n-\12)\kappa_0'+\omega_0.\end{array}\]
(2) Transmutation properties
\[\begin{array}{rl}
\p_z\Big({\red \cH}\big(\sigma,\kappa_n\big)+\omega_n\Big)&=
\Big({\red \cH}\big(\sigma,\kappa_{n+1}\big)
+\omega_{n+1}\big)\p_z,\\[3mm]
\Big(\sigma(z)\p_z+\kappa_{n+1}(z)\Big)
\Big({\red \cH}\big(\sigma,\kappa_{n+1}\big)+\omega_{n+1}\Big)&
=\Big({\red \cH}\big(\sigma,\kappa_n\big)+\omega_n\Big)
\Big(\sigma(z)\p_z+\kappa_{n+1}(z)\Big)
.\end{array}\]
(3) Recurrence relations
\[\begin{array}{rll}&
\Big({\red \cH}\big(\sigma,\kappa_n\big)+\omega_n\Big)F=0\
&\Rightarrow\qquad
\Big({\red \cH}\big(\sigma,\kappa_{n+1}\big)+
\omega_{n+1}\Big)\partial_z F=0,\\[3mm]
&
\Big({\red \cH}\big(\sigma,\kappa_{n+1}\big)
+\omega_{n+1}\Big)F=0\ &\Rightarrow\qquad
\Big({\red \cH}\big(\sigma,\kappa_n\big)+\omega_n\Big) \Big(\sigma(z)\p_z+\kappa_{n+1}(z)\Big)F=0
.\end{array}\]\label{recr}
\eet


The argument used in the proof of the implication
(1)$\Rightarrow$(2) is typical for the so-called {\em supersymmetric
  quantum mechanics} and is described in the following lemma:

\bel
Suppose that for $n\in\cc$,  $\A_+^{n+\12}$, $\A_-^{n-\12}$, ${\red \cH}^n$ are  operators
and $\eta_{n\pm\12}$ are numbers satisfying
\[\begin{array}{rl}
{\red \cH}^n&=\A_+^{n+\12}\A_-^{n-\12}+\eta_{n-\12},\\[3mm]
&=\A_-^{n-\12}\A_+^{n+\12}+\eta_{n+\12}.\end{array}\]
Then 
\[\begin{array}{l}
\A_-^{n-\12}{\red \cH}^n={\red \cH}^{n-1}\A_-^{n-\12},\\[3mm]
\A_+^{n+\12}{\red \cH}^n={\red \cH}^{n+1}\A_+^{n+\12}.
\end{array}\]
\label{rer}\eel

\noindent
{\bf Proof of Theorem \ref{recr}}
  (1) follows by a direct computation.
 (1)  implies (2) by Lemma \ref{rer} if we set
\beq
\begin{array}{l}
\A_+^{n+\12}=\p_z,\ \ \ \ 
\A_-^{n-\12}=\sigma(z)\p_z+\kappa_n(z),\ \ \ \ 
\eta_{n+\12}=n(n+1)\frac{\sigma''}{2}+(n+\12)\kappa'+\omega,\\[3mm]
  {\red \cH}^n:={\red \cH}(\sigma,\kappa_n)+\omega_n.
\end{array}\label{anni}\eeq
 (2) easily 
implies (3). 
\qed

\subsection{Ladders of solutions}
\label{Ladders of solutions}

Let $\kappa_n,\omega_n$ be as in
\eqref{kappan} and  \eqref{lambdan}.
 {\red Let $\rho_0$ solve 
  \eqref{e1a} for  $\sigma_0,\kappa_0$, that is, 
  \beq\big(\sigma_0(z)\p_z-\kappa_0(z)\big)\rho_0(z)=0.\label{e1b.}\eeq}
Suppose we have a solution
\beq
\Big({\red \cH}(\sigma,\kappa_0)+\omega_0\Big)f=0.\eeq
Then from $f$ we can construct  a two-sided ladder of
solutions. More precisely, for any $n\in\nn_0$ 
\begin{align}\label{ladder1}
 {\red \cH}(\sigma,\kappa_n)+\omega_n\quad\text{annihilates}\quad &\p_z^n f,\\
    {\red \cH}(\sigma,\kappa_{-n})+\omega_{-n}\quad\text{annihilates}\quad&\sigma^n\rho_0^{-1}\p_z^n\rho_0
                                                                     f
.\label{ladder2}
\end{align}

To see \eqref{ladder2}, we note that
\begin{align}
  \sigma(z)\p_z+\kappa_{-j}(z)&=\sigma(z)^{j+1}\rho_{\red
                                0}(z)^{-1}\p_z^j\sigma(z)^{-j}\rho_{\red
                                0}(z),\label{anni1-}\\
  \text{hence }\big(\sigma(z)\p_z+\kappa_{-n}(z)\big)\cdots
  \big(\sigma(z)\p_z+\kappa(z)\big)&=\sigma(z)^n\rho_0(z)^{-1}\p_z^n\rho_0(z).\label{anni2-} \end{align}

Consider now special cases $\omega_0=\mp\frac{\kappa_0'}{2}$:
\begin{align}
  {\red \cH}(\sigma,\kappa_0)-\frac{\kappa_0'}{2}&=\sigma(z)\p_z^2+\big(\sigma'(z)+\kappa_0(z)\big)\p_z,\\
  {\red \cH}(\sigma,\kappa_0)+\frac{\kappa_0'}{2}&=\sigma(z)\p_z^2+\big(\sigma'(z)+\kappa_0(z)\big)\p_z+\kappa_0'.
\end{align}
They have elementary solutions:
\begin{align}
  {\red \cH}(\sigma,\kappa_0)-\frac{\kappa_0'}{2}\qquad \text{ annihilates }&1,\label{anni3.}\\
  {\red \cH}(\sigma,\kappa_0)+\frac{\kappa_0'}{2} \qquad\text{ annihilates }&\rho_0^{-1}.\label{anni3}
                                      \end{align}

                                      \eqref{anni3.} is obvious. 
To see \eqref{anni3} we differentiate
\[(\sigma(z)\p_z+\kappa_0(z))\rho_0^{-1}(z)=0,\]
obtaining
\beq
0=(\sigma(z)\p_z^2+(\sigma'(z)+\kappa_0(z))\p_z+\kappa_0')\rho_0^{-1}(z)\\[3mm]
\eeq
Alternatively, we can derive \eqref{anni3} from \eqref{anni3.} by the
{\red basic}
symmetry \eqref{syme0}. {\red Indeed, 
\eqref{anni3.} and \eqref{syme0} imply that
$\cH(\sigma,-\kappa_0)-\frac{\kappa_0'}{2}$ annihilates
$\rho_0$. Then we switch the sign of $\kappa_0$, which corresponds
to replacing $\rho_0$ with $\rho_0^{-1}$.}

The special solutions \eqref{anni3.} and
\eqref{anni3} lead to a pair of one-sided ladders:
 for any $n\in\nn_0$
 \begin{align}\label{ladder2a}
  {\red \cH}(\sigma,-n\sigma'+\kappa_0)-(n+\tfrac12)\kappa_0'+n^2\frac{\sigma''}{2}
  \quad\text{annihilates}\quad&\sigma^n\rho_0^{-1}\p_z^n\rho_0
\\
   \label{ladder1a}
  {\red \cH}(\sigma,n\sigma'+\kappa_0)+(n+\tfrac12)\kappa_0'+n^2\frac{\sigma''}{2}\quad\text{annihilates}\quad &\p_z^n \rho_0^{-1}.
\end{align}
\eqref{ladder2a} consists of polynomials, which we will analyze in
more detail in Section \ref{Hypergeometric  polynomials}.
{\red The ladder \eqref{ladder1a} consists of functions of the form
$\rho_0^{-1}\sigma^{-n}P_n$, where
$P_n:=\sigma^n\rho_0\partial_z^n\rho_0^{-1}$ are polynomials.
}

\subsection{Chebyshev ladder}
\label{Chebyshev functions}

Let us fix $\sigma$, as usual a polynomial of degree  $\leq2$,  and
$\omega\in\cc$.
The equations given by the following two operators can be easily solved
in elementary functions:
\begin{align}\label{pwr1}
  {\red \cH}\Big(\sigma,-\frac{\sigma'}{2}\Big)+\omega+\frac{\sigma''}{4}&=\sqrt{\sigma(z)}\p_z\sqrt{\sigma(z)}\p_z+\omega\\
                                                                   &=\sigma(z)\p_z^2+\frac{\sigma'(z)}{2}\p_z+\omega,\\\label{pwr2}
  {\red \cH}\Big(\sigma,\frac{\sigma'}{2}\Big)+\omega+\frac{\sigma''}{4}&=\p_z\sqrt{\sigma(z)}\p_z\sqrt{\sigma(z)}+\omega\\
                                                                   &=\sigma(z)\p_z^2+\frac{3\sigma'(z)}{2}\p_z+\frac{\sigma''}{2}+\omega.
  \end{align}
  In fact, set
  \beq
  y(z)=\int_0^z\frac{\d x}{\sqrt{\sigma(x)}},\eeq
  which solves the equation $\frac{\d y}{\d
    z}=\frac{1}{\sqrt{\sigma(z)}}$.
  Then
  \begin{align}
    {\red \cH}\Big(\sigma,-\frac{\sigma'}{2}\Big)+\omega+\frac{\sigma''}{4}&\qquad
                                                                      \text{annihilates}\qquad\Big(A\sin\big(\omega
                                                                       y(z)\big)+B\cos\big(\omega
                                                                       y(z)\big)\Big),\\
        {\red \cH}\Big(\sigma,\frac{\sigma'}{2}\Big)+\omega+\frac{\sigma''}{4}&\qquad
                                                                      \text{annihilates}\qquad\frac{1}{\sqrt{\sigma(z)}}\Big(A\sin\big(\omega y(z)\big)+B\cos\big(\omega y(z)\big)\Big).
\end{align}

We can embed the parameters of \eqref{pwr1} and \eqref{pwr2} into a
single ladder as follows.
We set $\kappa_0=0$ and $\omega_0=\omega{\red +}\frac{\sigma''}{8}$, so that
\beq \kappa_m=m\sigma',\quad
\omega_m=\omega+\Big(m^2{\red+}\frac14
\Big)\frac{\sigma''}{2}.\eeq
then
\begin{align}
  \eqref{pwr1}&={\red \cH}\big(\sigma,\kappa_{-\frac12}\big)+\omega_{-\frac12},\\
  \eqref{pwr2}&={\red \cH}\big(\sigma,\kappa_{\frac12}\big)+\omega_{\frac12}.
                \end{align}
                By using the recurrence relations \eqref{ladder1} and \eqref{ladder2}
                we obtain elementary
                solutions for
                {\red \beq
                 \cH(\sigma,\kappa_m)+\omega_m
=\sigma(z)\partial_z^2+(m+1)\sigma'(z)\partial_z+\Big(m+\frac12\Big)^2\frac{\sigma''}{2}+\omega
\eeq
with $m\in\zz+\frac12$.

\begin{remark} Chebyshev ladders are very special, since $\kappa$ is
determined by $\sigma$. The only
two nontrivial hypergeometric types where they appear are 
${}_2F_1$  and $F_1$. However, they are quite important. One can note
that  Bessel
functions of half-integer parameters are elementary functions because
they correspond to the Chebyshev ladder of the $F_1$ type.\end{remark}
}

\init
\section{Singular point}
\label{Singular point}

Throughout this section we assume that 
\beq \sigma(0)=0.\label{singul}\eeq
We will discuss solutions of hypergeometric class equations
given by power series around $0$. We will also describe some
symmetries which exist if  \eqref{singul} holds.

\subsection{Solutions around a singular point}
\label{s.7}

\eqref{singul} usually means that $0$ is a singular point of the equation
 (\ref{e1}).
 In fact, rewriting the equation  (\ref{e1}) as
 \beq\frac{1}{\sigma(z)}\big({\red \sigma(z)\partial_z^2+\tau(z)\partial_z}+\omega\big)f(z)=\Big(\partial_z^2+\frac{\tau(z)}{\sigma(z)}\partial_z+\frac{\omega}{\sigma(z)}\Big)f(x)=0\eeq 
   we see that $\frac{\tau(z)}{\sigma(z)}$ or
   $\frac{\omega}{\sigma(z)}$ will usually have a singularity at $0$.
   Straightforward calculations 
                                (known under the name of the {\em
 Frobenius method} \cite{WW})
lead then to the following result
\bet There exists a unique formal power series $F(z)$
 solving 
\beq\Big({\red \cH}(\sigma,\kappa)+\omega\Big)F(z)=0,\ \ \ F(0)=1.\eeq
It is equal to
\beq F(\sigma,\kappa,\omega;z)=\sum_{n=0}^\infty
\frac{\Pi_{j=0}^{n-1}(\omega+(j+\12)\kappa'+j(j+1)\frac{\sigma''}{2})}
{\Pi_{j=0}^{n-1}(\kappa(0)+(j+1)\sigma'(0)) n!}(-z)^n
\label{qw1}\eeq
\eet

If $\sigma'(0)=0$ and $\sigma''\neq0$, then  $\sigma(z)$ has a double
root at zero and  then the series
(\ref{qw1})
 {\red either terminates and defines a polynomial, or
is divergent. However, it is always} asymptotic to
 one of the solutions of (\ref{e1}) {\red defined on
   $\cc\backslash[0,\infty[$ (see Section \ref{The  ${}_2F_0$
   equation}  and Appendix \ref{appendix}}).
In all other cases 
 the series
 (\ref{qw1}) {\red has a nonzero radius of convergence}. The  function
 $F(\sigma,\kappa,\omega;z)$ given by \eqref{qw1}
 will be called the {\em unified hypergeometric function}.
 It depends meromorphically on $\sigma'(0),\sigma'',\kappa,\omega$.

Recall that $\kappa_n$ and $\omega_n$ were defined in \eqref{kappan}
and \eqref{lambdan}. The basic pair of recurrence relations 
is as follows:
\begin{align}\label{pair1}
\p_z F(\sigma,\kappa_n,\omega_n;z)&
=-\frac{\omega_n+\frac12\kappa_n'}{\kappa_{n+1}(0)}
 F(\sigma,\kappa_{n+1},\omega_{n+1};z),\\
\left(\sigma(z)\p_z+\kappa_{n+1}(z)\right)
 F(\sigma,\kappa_{n+1},\omega_{n+1};z)&=
\kappa_{n+1}(0)F(\sigma,\kappa_n,\omega_n;z).\label{pair2}
\end{align}

\subsection{Olver's normalization}
\label{s.7a}

Assume now $\red\sigma'(0)\neq0$. We set $m:=\frac{\kappa(0)}{\sigma'(0)}$.
Without sacrificing the generality we can suppose additionally that $\sigma'(0)=1$.
Thus\beq\sigma(z)=z+\frac{\sigma''}{2}z^2,\qquad \kappa(z)=m+\kappa' z.\label{olver.1}\eeq
$0$ is a regular singular point and $0$, $-m$ are its indices, see
\cite{WW} or Appendix \ref{Class}.

In this case often instead of the function (\ref{qw1}) it is more convenient
to use the function
\beq
{\mathbf F}(\sigma,\kappa,\omega;z):=\frac{F(\sigma,\kappa,\omega;z)}{\Gamma(1+m)}=
\sum_{n=0}^\infty
\frac{\Pi_{j=0}^{n-1}(\omega+(j+\12)\kappa'+j(j+1)\frac{\sigma''}{2})}
{\Gamma(m+n+1)n!}(-z)^n.\label{olver}\eeq
Note that ${\mathbf F}$ is holomorphic in $\sigma'',\kappa,\omega$.
Dividing by $\Gamma(1+m)$, as in \eqref{olver}, is sometimes called
{\em Olver's normalization}, {\red since it was made popular by Frank Olver's
textbook \cite{Ol2}.}
Here is the basic pair of recurrence relations for the {\em Olver
  normalized unified hypergeometric function}:
\begin{align}\label{pair3}
\p_z{\mathbf  F}(\sigma,\kappa_n,\omega_n;z)&
=-\Big(\omega_n+\frac12\kappa_n'\Big)
 {\mathbf F}(\sigma,\kappa_{n+1},\omega_{n+1};z),\\
\left(\sigma(z)\p_z+\kappa_{n+1}(z)\right)
    {\mathbf F}(\sigma,\kappa_{n+1},\omega_{n+1};z)&=
                                                                         {\mathbf F}(\sigma,\kappa_n,\omega_n;z).\label{pair4}
\end{align}

\subsection{Power symmetry}
\label{sec.sec}

We keep the assumptions of Subsection \ref{s.7a}.
Let us introduce 
a certain transformation of the  parameters.  We keep $\sigma$ the same, and the transformed $\kappa$ and
$\omega$ are 
\beq
\kappa^\circ(z):=-m+(\kappa'-m\sigma'')z,\quad\omega^\circ:=\omega-m\kappa'+m^2\frac{\sigma''}{2}.\label{olver.2}\eeq
Consequently, $m^\circ=-m$.

\bet
\beq
z^m\Big({\red \cH}(\sigma,\kappa)+\omega\Big)z^{-m}={\red \cH}(\sigma,\kappa^\circ)+\omega^\circ.\label{powers}\eeq
\eet

\proof
Using
\beq z^m\p_z z^{-m}=\p_z-\frac{m}{z},\qquad
z^m\p_z^2 z^{-m}=\p_z^2-\frac{2m}{z}\p_z+\frac{m(m+1)}{z^2}
\eeq
we compute
\begin{align}
  &z^m\Big({\red \cH}(\sigma,\kappa)+\omega\Big)z^{-m}\\
  =&z^m\Big(\Big(z+\frac{\sigma''}{2}z^2\Big)\p_z^2+\big(1+m+(\sigma''+\kappa')z\big)\p_z+\frac{\kappa'}{2}+\omega\Big)z^{-m}\\
  =&\Big(z+\frac{\sigma''}{2}z^2\Big)\p_z^2+\Big(1-m+\big((1-m)\sigma''+\kappa'\big)z\Big)\p_z+\frac{\sigma''}{2}m(m-1)-\kappa'(m-\tfrac12)+\omega\\=&
  {\red \cH}(\sigma,\kappa^\circ)+\omega^\circ.
  \end{align}

  Note that $\kappa^{\circ\circ}=\kappa$ and 
$\omega^{\circ\circ}=\omega$. Thus we obtain a $\zz_2$ symmetry
of {\red the hypergeometric class.  It is} different from the basic symmetry
of Thm \ref{syme}.

As a consequence of \eqref{powers},  both $\mathbf{F}(\sigma,\kappa,\omega;z)$ and $z^{-m}{\mathbf
  F}(\sigma,\kappa^\circ,\omega^\circ;z)$ are annihilated by $ 
{\red \cH}(\sigma,\kappa)+\omega$. If $m\not\in\zz$ they form a basis of
solutions of the corresponding equation. The situation for $m\in\zz$
will be discussed in Section \ref{Degenerate case}
about the degenerate case.

  \subsection{Inversion symmetry}

  We still assume $\sigma(0)=0$. 
  Besides, we suppose that
  \beq \label{poset1}\sigma''\neq0\quad\text{ or }\quad \kappa'\neq0.\eeq
  Suppose that $\zeta\in\cc$ solves the equation
  \beq
  \frac{\sigma''}{2}\zeta^2+(\sigma''+
  \kappa')\zeta+\frac{\kappa}{2}+\omega=0.\label{poset}\eeq
  Note that \eqref{poset1} guarantees that \eqref{poset} has a
  solution. Set
  \begin{align}
    \sigma^\triangle&:=\frac{\sigma''}{2}w-\sigma'(0)w^2,\\
    \kappa^\triangle&=\big(\kappa(0)+2(\zeta+1)\sigma'(0)\big)w-\sigma''(1+\zeta)-\kappa',\\
    \omega^\triangle&:=-\sigma'(0)(1+\zeta)^2-\kappa(0)(\zeta+\tfrac12).
    \end{align}
    \bet
      Consider the substitution $w=-\frac{1}{z}$. Then we have 
  \begin{align}
    &-z^{-\zeta+1}\Big({\red \cH}(\sigma,\kappa;z,\p_z)+\omega\Big)z^{\zeta}
      ={\red \cH}(\sigma^\triangle,\kappa^\triangle;w,\p_w)+\omega^\triangle\end{align}
    Hence
    \beq\label{inverso}
    {\red \cH}(\sigma,\kappa)+\omega\quad\text{ annihilates }\quad z^{-\zeta}\mathbf{F}(\sigma^\triangle,\kappa^\triangle,\omega^\triangle;-z^{-1}).\eeq
\label{inversion}\eet

\proof
    Indeed, using $\p_z=w^2\p_w$ and $\p_z^2=w^4\p_w^2+2w^3\p_w$, we obtain
  \begin{align}
    &    -z\Big({\red \cH}(\sigma,\kappa;z,\p_z)+\omega\Big)\\
    =&\sigma(-w^{-1})w^3\p_w^2+\Big(2\sigma(-w^{-1})w^2+\sigma'(-w^{-1})w+\kappa(-w^{-1})w\Big)\p_w+\Big(\frac{\kappa'}{2}+\omega\Big)w^{-1}.
       \label{porru}\end{align}
     Then we use $w^\zeta\p_ww^{-\zeta}=\p_w-\frac{\zeta}{w}$ and
     $w^\zeta\p_w^2w^{-\zeta}=\p_w^2-\frac{2\zeta}{w}\p_w+\frac{(\zeta+\zeta^2)}{w^2}$
     to obtain
     \begin{align}&
                    w^\zeta\eqref{porru}w^{-\zeta}\\
       =&\Big(\frac{\sigma''}{2}w-\sigma'(0)w^2\Big)\p_w^2+\Big(\big(\kappa(0)+2\zeta\sigma'(0)\big)w-\frac{\sigma''}{2}-\zeta\sigma''-\kappa'\Big)\p_w\\
                  &-\sigma'(0)(\zeta+\zeta^2)-\kappa(0)\zeta+\Big(\frac{\kappa'}{2}+\omega+\frac{\sigma''}{2}(2\zeta+\zeta^2)+\kappa'\zeta\Big)w^{-1}\\
       =&\sigma^\triangle(w)\p_w^2+\big({\sigma^\triangle}'(w)+\kappa^\triangle(w)\big)\p_w+\frac{{\kappa^\triangle}'}{2}+\omega^\triangle.
  \end{align}
\qed

Note that
$\sigma^{\triangle \triangle}=\sigma$, $\kappa^{\triangle
  \triangle}=\kappa$, $\omega^{\triangle \triangle}=\omega$. Hence the
inversion symmetry  generates a group isomorphic to $\zz_2$ {\red
  acting on the hypergeometric class}.

\begin{remark} It is easy to see that the substitution $z\mapsto \frac1w$
also  generates a symmetry of hypergeometric class equations. This follows
  immediately from Theorem \ref{inversion} and the trivial fact that
  $z\mapsto -z$ leads to a symmetry {\red of the hypergeometric class}
  as well. \end{remark}

\init
\section{Integral representations}
\label{Integral representations}

Hypergeometric class functions possess useful integral
representations.  Most of them have the form of an Euler 
transform of an elementary function. In Subsection \ref{s.9} we will
show how to derive these 
representations in a unified way. Note that typically for a given
equation   one can choose various contours of integration, obtaining various
solutions.

Some important integral represenations have a different form---they
can be viewed as Laplace transforms of certain elementary functions.
They will be described in Subsection \ref{s.10}.

In this section we will often deal with
multivalued functions $s\mapsto f(s)$ defined
on a certain  Riemann surface.
These functions are analytically continued along a certain curve $\gamma$
contained in this Riemann surface.
We will use the notation
\beq f(s)\Big|_{s_0}^{s_1}=f(s_1)-f(s_0),\eeq
where $s_1$ and $s_2$ are the endpoints of the curve $\gamma$.

As usual, we fix $\sigma,\kappa_0$, and for $n\in\cc$, 
as in \eqref{kappan}, we set \[\kappa_n(z):=n\sigma'(z)+\kappa_0(z).\]
In this 
section  $\omega_0=\frac{\kappa_0'}{2}$, so that according to 
\eqref{lambdan}, 
 \beq 
\omega_n:=n^2\frac{\sigma''}{2}+(n+\tfrac12)\kappa_0'.\label{keep}\eeq

\subsection{Euler transforms}
\label{s.9}

As usual, we assume that 
\beq 
(\sigma(z)\p_z-\kappa_0(z))\rho_0(z)=0.\label{vanish4}
\eeq

\bet Let $n\in\cc$.
Suppose that the curve $[0,1]\ni \tau\overset{\gamma}\mapsto s(\tau)$ satisfies
\beq
\sigma(s)(s-z)^{-n-2}\rho_0^{-1}(s)\Big|_{\s(0)}^{s(1)}=0
\label{scon}\eeq
and 
\beq
f_n(z):=\int_\gamma
(s-z)^{-n-1}\rho_0^{-1}(s)\d s
\label{q2}\eeq
is well-defined. Then $f_n$ is annihillated by
\begin{align}\label{vanish3}
&{\red \cH}(\sigma,\kappa_n)+ \omega_n\\
=&\sigma(z)\p_z^2+(\kappa_0(z)+(n+1)\sigma'(z))\p_z+
\kappa_0'(n+1)+\frac{\sigma''}{2}n(n+1)
.\end{align}
Besides, if $f_{n+1}$ is also well-defined, then
\begin{align}\label{recu1}
\p_zf_n(z)&=(n+1)f_{n+1}(z),\\\label{recu2}
\left(\sigma(z)\p_z+\kappa_{n+1}(z)\right)f_{n+1}(z)&
=-\left(\kappa_0'+\frac{\sigma''}{2}n\right)f_{n}(z).
\end{align}
\label{in1}\eet

\proof For simplicity, we will assume that both $f_n$ and $f_{n+1}$
are well defined. \eqref{recu1} is obvious. Let us prove
\eqref{recu2}. Using the Taylor expansion of the polynomials
$\sigma(z)$ and $\kappa(z)$ we obtain
\begin{align}
  &\Big(\sigma(z)\p_z+\kappa_0(z)+(n+1)\sigma'(z)\Big)f_{n+1}(z)\\
  =&(n+2) \int_\gamma\Big(\sigma(s)+\sigma'(s)(z-s)+\frac{\sigma''}{2}(z-s)^2\Big)
     (s-z)^{-n-3}\rho_0(s)^{-1}\d s\\
  &+\int_\gamma\Big(\kappa_0(s)+(z-s)\kappa'+(n+1)\sigma'(s)+(n+1)(z-s)\sigma''\Big)(s-z)^{-n-2}\rho_0(s)^{-1}\d
    s\\
  =&-\Big(n\frac{\sigma''}{2}+\kappa_0'\Big)\int_\gamma(s-z)^{-n-1}\rho_0(s)^{-1}\d
     s\\
  &+\int_\gamma\Big((n+2)\sigma(s)(s-z)^{-n-3}-\sigma'(s)(s-z)^{-n-2}+\kappa_0(s)(s-z)^{-n-2}\Big)\rho_0(s)^{-1}\d
    s\\
  =&-\Big(n\frac{\sigma''}{2}+\kappa_0'\Big)f_n(z)\\
  &+\int_\gamma\Big(\p_s\sigma(s)(s-z)^{-n-2}\rho_0(s)^{-1}\Big)\d s\label{vanish1}\\
  &+\int_\gamma(s-z)^{-n-2}\Big(\sigma(s)\partial\rho_0(s)^{-1}+\kappa_0(s)\rho_0(s)^{-1}\Big)\d s.\label{vanish2}
  \end{align}
\eqref{vanish1} vanishes because of 
\eqref{scon}. \eqref{vanish2} is zero by \eqref{vanish4}. \qed

It is sometimes useful to consider differently normalized Euler
integrals:
\beq
\mathbf{f}_n(z):=\frac{f_n(z)}{\Gamma(n+1)}=\frac{1}{\Gamma(n+1)}\int_\gamma
(s-z)^{-n-1}\rho_0^{-1}(s)\d s
\label{q2a}\eeq
Then the recurrence relations \eqref{recu1} and \eqref{recu2} are modified:
\beq\begin{array}{rl}
\p_z\mathbf{f}_n(z)&=\mathbf{f}_{n+1}(z),\\[2mm]
\left(\sigma(z)\p_z+\kappa_{n+1}(z)\right)\mathbf{f}_{n+1}(z)&
=-(n+1)\left(\kappa_0'+\frac{\sigma''}{2}n\right)\mathbf{f}_{n}(z).
\end{array}\label{rec.1}\eeq

\subsection{Laplace integrals}
\label{s.10}

Suppose that 
$\sigma''=0$. We still assume \eqref{keep},
which now can be rewritten as
\beq \omega_n=(n+\tfrac12)\kappa_0'.\eeq
Let
\beq\Big(\sigma(-\p_s)s+\kappa_0(-\p_s)\Big) 
\delta_0(s)=0,\label{delt}\eeq 
which defines up to a coefficient an elementary function 
$\delta_0$.

\bet Assume that $[0,1]\ni \tau\overset{\gamma}\mapsto s(\tau)$ satisfies
\beq
(s^{n+2}\sigma'+s^{n+1}\kappa')\delta_0(s)\e^{sz}\Big|_{s(0)}^{s(1)}=0,\label{broke1}\eeq
and the following integral exists:
\beq 
g_n(z)=\int_\gamma\delta_0(s)s^{n}\e^{zs}\d s.\label{laplace}\eeq
Then $g_n$ is annihillated by
\begin{align}
  &{\red \cH}(\sigma,\kappa_n)+\omega_n\\
  =\,&\sigma(z)\p_z^2+(\sigma'(z)(n+1)+\kappa_0(z))\p_z+(n+1)\kappa_0'
.\label{broke}\end{align}
Besides, if also $g_{n+1}$ is well defined, then
\begin{align}\label{reccu1}
\p_zg_n(z)&=g_{n+1},\\\label{reccu2}
(\sigma(z)\p_z+\kappa_{n+1}(z))g_{n+1}(z)&=
-(n+1)\kappa_0'g_{n}(z).
\end{align}
\label{in1-}\eet

\proof
For simplicity, we will assume that both $g_n$ and $g_{n+1}$ are well
defined, and we will prove the recurrence relations, which directly imply 
\eqref{broke}.

\eqref{reccu1} is obvious. Let us show \eqref{reccu2}:
\begin{align}
&  (\sigma(z)\p_z+(n+1)\sigma'+\kappa_0(z)) g_{n+1}\\
  =&\int_\gamma\Big(\sigma'\delta_0(s)s^{n+2}z+
     \sigma(0)\delta_0(s)s^{n+2}+(n+1)\sigma'\delta_0(s)s^{n+1}\Big)\e^{zs}\d 
     s\\&
  +\int_\gamma\Big(\kappa_0'\delta_0(s)s^{n+1}z+
     \kappa_0(0)\delta_0(s)s^{n+1}\Big)\e^{zs}\d 
  s\\
  =&\int_\gamma\d
     s\partial_s\Big(\sigma'\delta_0(s)s^{n+2}\e^{zs}+\kappa_0'\delta_0(s)s^{n+2}\e^{zs}\Big)\label{van1}\\
  &+\int_\gamma\Big(-\sigma's\delta_0'(s)-\sigma'\delta_0(s)-\kappa_0'\delta_0'(s)+\sigma(0)s\delta_0(s)+\kappa_0(0)\delta_0(s)\Big)s^{n+1}\e^{zs}\d
    s\label{van2}\\
  &-(n+1)\kappa_0'\int\delta_0(s)s^n\e^{zs}\d s.
  \end{align}
Now \eqref{van1} vanishes because of \eqref{broke1} and \eqref{van2}
due to \eqref{delt}. \qed

Sometimes it will be more convenient to present the integral represention
\eqref{laplace} in a different form, with the variable $s$ replaced
with $t=\frac1s$:
\beq 
g_n(z)=\int_{\tilde\gamma}\delta_0(t^{-1})t^{-n+2}\e^{\frac{z}{t}}\d t.\label{laplace1}\eeq
where the contour $\tilde\gamma$ is obtaned by the change of variable
and  reversing the orientation.
\init
\section{Applications case by case}
\label{Applications case by case}

{\red
Hypergeometric class operators can be divided into several types. By an affine transformation $z\mapsto az+b$ and
division by a constant, an
operator of each type  can be reduced to its normal form. There are 5
nontrivial types, with normal forms listed in the introduction.
For instance, if $\sigma$ has two distinct roots, then the equation
belongs to the ${}_2F_1$ type, if it has a double root, it belongs to
the ${}_2F_0$ type, etc.
The hypergometric class contains also 4 ``trivial types'', which can be solved in an
elementary way. All these 9 types are listed in the table in Appendix
\ref{Class}.

Strictly speaking, this table is devoted to types of the Riemann
class, which is larger than the hypergeometric class. Therefore, this
table contains 10 types. It includes one additional type: the Airy
operator, which cannot be reduced to the hypergeometric
class. However, all the 9 other types of the Riemann class are represented
in the hypergeometric class.

In
this section
we describe
the properties of all five nontrivial types of the hypergeometric class
 that are direct consequences of the unified theory,
discussed in  the previous sections.}
We try to follow the same pattern. Here is the list of items that we
will give for each hypergeometric  type, if available:
\begin{itemize}
\item Operator that generates the equation, that is
${\red \cH}(\sigma,\kappa)+\omega$,
  expressed in the traditional form.
\item Parameters $\sigma(z),\kappa(z),\omega$; see \eqref{war} and \eqref{coro}.
\item Weight $\rho(z)$, see \eqref{e1b.} or \eqref{e1b}.
\item Basic symmetry, see \eqref{asa}.
\item Power symmetry (if available), see \eqref{powers}.
\item Inversion symmetry (if available), see \eqref{inverso}.
  \item Theorem about integral representations, see Theorems \ref{in1}
    and \ref{in1-}.
    \item Standard solution annihillated by ${\red \cH}(\sigma,\kappa)+\omega$ (typically a
      special case of the unified hypergeometric function \eqref{qw1}, but
      for the Hermite equation we need to make an exception).
    \item Olver normalized standard solution (if available), see \eqref{olver}.
      \item Basic pair of recurrence relations for standard solutions,
        see \eqref{pair1}-\eqref{pair2} or        
        \eqref{pair3}-\eqref{pair4}.
        \item Integral representations of standard solution (obtained
          by an application of Theorems \ref{in1}
    or \ref{in1-}).
          \item Chebyshev solutions (if available), see Subsection \ref{Chebyshev functions}.
  \end{itemize}

  We will often use the {\em Pochhammer symbol}
  \beq(a)_j:=a(a+1)\cdots(a+j-1),\qquad a\in\cc,\quad j\in\nn_0.\eeq

\subsection{The ${}_2F_1$  equation}

\noindent
 {\bf ${}_2F_1$ operator:}
\begin{align}
  \mathcal{F}(a,b;c):=z(1-z)\p_z^2+\big(c-(a+b+1)z\big)\p_z-ab.\label{hy1-tra}
\end{align} 
\noindent  
{\bf Parameters:} 
\[
  \sigma(z)=z(1-z),\quad
  \kappa(z)=c-1-(a+b-1)z,\quad
  \omega=-\big(a-\tfrac12\big) \big(b-\tfrac12\big)-\tfrac14.\]
\noindent 
{\bf Weight:}
\[\rho(z)=z^{c-1}(z-1)^{a+b-c}.\]
                     {\bf Basic symmetry:}\[
 z^{c-1}(z-1)^{a+b-c}\F(a,b;c)z^{1-c}(z-1)^{-a-b+c}
=\F(1-b,1-a;2-c) .\]
{\bf Power symmetry:}
\[
 z^{c-1}\F(a,b;c)z^{1-c}=
\F(b+1-c,a+1-c;2-c) .\]
{\bf Inversion symmetry:}
\[
(-z)^{1+a}\F(a,b;c;z,\p_z)(-z)^{-a}=
\F(a,a-c+1;a-b+1;w,\p_w)
,\qquad w=z^{-1}.\]

  \bet[about integral representations]
Let $[0,1]\ni \tau\overset{\gamma}\mapsto t(\tau)$ satisfy 
\[t^{b-c+1}(1-t)^{c-a}(t-z)^{-b-1}\Big|_{\gamma(0)}^{\gamma(1)}=0.\]
Then
\beq 
\F(a,b;c) \quad\text{ annihilates }\quad
\int_\gamma t^{b-c}(1-t)^{c-a-1}(t-z)^{-b}\d t.
\label{f4}\eeq
\label{intr}\eet

\proof 
We check that 
\beq
\F(a,b;c) 
t^{b-c}(1-t)^{c-a-1}(t-z)^{-b}=-b 
\p_t t^{b-c+1}(1-t)^{c-a}(t-z)^{-b-1}
.\eeq
\qed\\
\noindent
{\bf ${}_2F_1$  function:}  For $|z|<1$ defined  by the
following power series, then
analytically extended:
\[F(a,b;c;z)=\sum_{j=0}^\infty
\frac{(a)_j(b)_j}{
(c)_j}\frac{z^j}{j!}, \quad |z|<1.\]
 {\bf Olver normalized ${}_2F_1$  function}:
\[ {\bf F}  (a,b;c;z):=\frac{F(a,b,c,z)}{\Gamma(c)}
=\sum_{j=0}^\infty
\frac{(a)_j(b)_j}{
\Gamma(c+j)}\frac{z^j}{j!}.\]

\noindent
{\bf Integral representation:}
\begin{eqnarray}\label{eqa1}
&&\int_1^\infty t^{b-c}(t-1)^{c-a-1}(t-z)^{-b}\d t\\
&=&
 \Gamma(a)\Gamma(c-a)
\mathbf{F}(a,b;c;z),\ \ \ \ \Re(c-a)>0,\ \Re a>0,\ \ \ z\not\in[1,\infty[.
\nonumber\end{eqnarray}

\noindent
{\bf Basic pair of recurrence relations:} 
\begin{align}\p_z \mathbf{F}(a,b;c;z)&=ab
                                       \mathbf{F}(a+1,b+1;c+1;z),\\
\Big( z(1-z)\p_z+\big(c-(a+b+1)z\big) \Big)\mathbf{F}(a+1,b+1;c+1;z)&=\mathbf{F}(a,b;c;z).
\end{align} 

\noindent
{\bf Chebyshev solutions} for $k=0,1,\dots$:
\begin{align}
{\bf F}\Big(1+ k+\lambda,1+k-\lambda;\frac32+k;\frac{1-w}{2}\Big)
 &= \frac{(1-w^2)^{\frac12+k}}{2\sqrt\pi(-2)^k}\partial_w^k
 \frac{\big(w+\ii\sqrt{1-w^2}\big)^\lambda 
     +\big(w-\ii\sqrt{1-w^2}\big)^\lambda}{\sqrt{1-w^2}},\\
{\bf F}\Big(- k+\lambda,-k-\lambda;\frac12-k;\frac{1-w}{2}\Big)
 = &\frac{2^{k}}{\ii\sqrt\pi(\lambda-k)_{2k+1}}
 \partial_w^k 
 \frac{\big(w+\ii\sqrt{1-w^2}\big)^\lambda 
          -\big(w-\ii\sqrt{1-w^2}\big)^\lambda}{\sqrt{1-w^2}} .
      \end{align}

\subsection{The ${}_1F_1$ equation}
\label{The ${}_1F_1$ equation}

\noindent {\bf
 ${}_1F_1$ operator}
\beq
\F(a;c):=z\p_z^2+(c-z)\p_z-a.
\label{f1c}\eeq
\noindent {\bf Parameters}
\begin{align*}
  \sigma(z)&=z,\quad
  \kappa(z)=c-1-z,\quad
  \omega=-a+\tfrac12.
\end{align*}
\noindent {\bf Weight:}
\[ \rho(z)=\e^{-z}z^{c-1}
.\]

\bet[about integral representations]\ben\item
Let $[0,1]\ni \tau\overset{\gamma}\mapsto t(\tau)$ satisfy
\[t^{a-c+1}\e^t(t-z)^{-a-1}\Big|_{t(0)}^{t(1)}=0.\]
Then
\beq
\F(a;c) \quad\text{ annihilates }\quad
\int_\gamma t^{a-c}\e^t(t-z)^{-a}\d t.\label{dad1}\eeq
\item
Let $[0,1]\ni \tau\overset{\gamma}\mapsto t(\tau)$ satisfy 
\[\e^{\frac{z}{t}}t^{-c}(1-t)^{c-a}\Big|_{t(0)}^{t(1)}=0.\]
Then
\beq
\F(a;c) \quad\text{ annihilates }\quad\int_\gamma\e^{\frac{z}{t}}t^{-c}(1-t)^{c-a-1}\d t
.
\label{dad}
\eeq
\een\label{dad4}\eet

\proof
\begin{align}
\F(a,c) \, t^{a-c}\e^t(t-z)^{-a}
&=-a\p_t
t^{a-c+1}\e^t(t-z)^{-a-1},\\
\F(a,c)\,\e^{\frac{z}{t}}t^{-c}(1-t)^{c-a-1}&=-\p_t
\e^{\frac{z}{t}}t^{-c}(1-t)^{c-a}.
\end{align}\qed


\noindent{\bf Basic symmetry:}\[
z^{c-1}\e^{-z}
\F(a;c;z,\p_z)z^{1-c}\e^z=
\F(1-a;2-c;w,\p_w),\quad w=-z.\]
{\bf Power symmetry:}\[
z^{c-1}
\F(a;c)z^{1-c}=
                \F(1+a-c;2-c).
                \]
 {\bf Inversion symmetry:}\[
                                        z^{a+1}\F(a;c;z,\p_z)z^{-a}=\cF(a,1+a-c;{\red-};w,\p_w),\qquad w=-z^{-1}.
\]
(Thus the ${}_1F_1$ equation is equivalent to the  ${}_2F_0$ 
equation.) \\
\noindent
{\bf ${}_1F_1$ function:}
\[F(a;c;z):=\sum_{j=0}^\infty
\frac{(a)_j}{
(c)_j}\frac{z^j}{j!}.\]
\noindent
 {\bf Olver normalized  ${}_1F_1$ function}
\[ {\bf F}  (a;c;z):=\frac{F(a;c;z)}{\Gamma(c)}=
\sum_{j=0}^\infty
\frac{(a)_j}{
\Gamma(c+j)}\frac{z^j}{j!}.\]

\noindent{\bf Basic pair of recurrence relations:}
\begin{align}
  \p_z\mathbf{F}(a;c;z)&=a\mathbf{F}(a+1;c+1;z),\\
  \big( \p_z+c-z)\big)\mathbf{F}(a+1;c+1;z)&=\mathbf{F}(a;c;z).
                                             \end{align}
{\bf Integral representations:} for all parameters
\begin{eqnarray*}
\frac{1}{2\pi \i}\int\limits
_{]-\infty,(0,z)^+,-\infty[} t^{a-c}\e^t(t-z)^{-a}\d t
&=& {\bf F}  (a;c;z);\end{eqnarray*}
for $\Re a>0,\ \Re (c-a)>0$
\begin{eqnarray*}\int\limits_{[1,+\infty[}\e^{\frac{z}{t}}t^{-c}(t-1)^{c-a-1}\d t
&=& \Gamma(a)\Gamma(c-a)\mathbf{F}(a;c;z).\end{eqnarray*}

\subsection{The  ${}_2F_0$ equation}
\label{The  ${}_2F_0$ equation}

\noindent{\bf The  ${}_2F_0$ operator}:
\beq {\F}(a,b;-):=z^2\partial_z^2 +( -1+(1+a+b)z)\partial_z +ab, \label{g8}\eeq
\noindent {\bf Parameters:} 
\[
\sigma(z)=z^2,\quad
\kappa(z)=-1+(a+b-1)z,\quad
  \omega=\big(a-\tfrac12\big) \big(b-\tfrac12\big)+\tfrac14.
\]
\noindent {\bf Weight:} 
\[\rho(z)=z^{-1+a+b}\e^{\frac1z}.
\]
{\bf Basic symmetry:}\[
  z^{-1+a+b}\e^{\frac1z}\F(a,b;-;z,\p_z)z^{1-a-b}\e^{-\frac1z}
  =\cF(1-b,1-a;-;w,\p_w),\quad w=-z.\]
 {\bf Inversion symmetry:}\[
z^{a+1}\F(a,b;-;z,\p_z)z^{-a}
=-\F(a;1+a-b;w,\p_w),\quad w=-z^{-1}.\]
{\red(Thus the ${}_2F_0$ equation is equivalent to the  ${}_1F_1$ 
equation.) }

\bet[about integral representations]
Let $[0,1]\ni \tau\overset{\gamma}\mapsto t(\tau)$ satisfy 
\[\e^{-\frac{1}{t}}t^{b-a-1}(t-z)^{-b-1}
\Big|_{t(0)}^{t(1)}=0.\]
Then
\beq\F(a,b;-) \quad\text{ annihilates }\quad\int_\gamma\e^{-\frac{1}{t}}t^{b-a-1}(t-z)^{-b}\d t.
\label{dad2}\eeq
A second integral representation is obtained if we interchange $a$ and $b$.  \eet

\proof
We check that 
\beq
\F(a,b;-)\e^{-\frac{1}{t}}t^{b-a-1}(t-z)^{-b}=-b\p_t\e^{-\frac{1}{t}}t^{b-a-1}(t-z)^{-b-1}.\eeq
\qed\\\noindent
{\bf ${}_2F_0$ function:} It is defined
for $z\in\cc\backslash[0,+\infty[$,
\[F(a,b;-;z):=\lim_{c\to\infty}F(a,b;c;cz),\]
where $|\arg c|>\epsilon$, $\epsilon>0$.
It extends to an analytic function on the universal cover of
$\cc\backslash\{0\}$ 
with a branch point of an infinite order at 0. It is annihilated by $\cF(a,b;-)$.\\
{\bf Asymptotic expansion}:
\[
F(a,b;-;z)\sim\sum_{j=0}^\infty\frac{(a)_j(b)_j}{j!}z^j,
\ \red |\arg z|>\epsilon.
\]
\noindent{\bf Basic pair of recurrence relations:}
\begin{align}
  \p_zF(a,b;-;z)&=abF(a+1,b+1;-;z)\\
  \Big( z^2\p_z+\big(-1+(1+a+b)z)\big)\Big)F(a+1,b+1;-;z)&=F(a,b;-;z).
                                             \end{align}
{\bf Integral representation} for $\Re a>0$:
\[\int_0^\infty 
\e^{-\frac{1}{t}}t^{b-a-1}(t-z)^{-b}\d t 
= \Gamma(a)F   (a,b;-;z), \ \ z\not\in[0,\infty[.
\]

{\red
The function ${}_2F_0$ is in our opinion insufficiently
known. Therefore, we devote Appendix \ref{appendix} to a derivation of
some of its properties.

  Note that the
equivalence of ${}_2F_0$ and ${}_1F_1$ equations is helpful in deriving
properties of Subsection \ref{The ${}_1F_1$ equation}
from Subsection \ref{The  ${}_2F_0$ equation}. For instance, the
integral representation \eqref{dad2} can be deduced from the integral
representation \eqref{dad} and the inversion symmetry. However, the
relationship between some of the properties is not so
straightforward.
For instance, the basic recurrence relations for
${}_2F_0$ equations  do not follow directly from the basic recurrence
relations for the  ${}_1F_1$ equations. Besides
${}_2F_0$ equation
does not possess the power symmetry,
  even though the  ${}_1F_1$ equation has it. 
  Instead, the ${}_2F_0$ equation is symmetric with respect to the
  parameter interchange
  $a\leftrightarrow b$.}

\subsection{The ${}_0F_1$ equation}

\noindent{\bf ${}_0F_1$  operator:}
\[\F(c;z,\p_z):=z\p_z^2+c\p_z-1.\]
\noindent {\bf Parameters:} 
\[
  \sigma(z)=z,\quad
  \kappa(z)=c-1,\quad
  \omega=-1. \]
\noindent {\bf Weight:} 
\[  \rho(z)=z^{c-1}.
           \]
\noindent{\bf\red Basic symmetry $=$ power symmetry:}  \begin{align}             z^{c-1}
\F(c)z^{1-c}&=
\F(2-c). 
\label{tozs1}\end{align}

\bet[about integral representations]\label{schl}
Suppose that $[0,1]\ni t\mapsto\gamma(t)$ satisfies
\[\e^t\e^{\frac{z}{t}}t^{-c}\Big|_{\gamma(0)}^{\gamma(1)}=0.\]
Then
\beq\F(c;z,\p_z) \quad\text{ annihilates }\quad
\int_\gamma\e^t\e^{\frac{z}{t}}t^{-c}\d t.\label{dad5}\eeq
\eet

\proof
We check that 
\beq\F(c)\,
\e^t\e^{\frac{z}{t}}t^{-c}=-
\p_t\e^t\e^{\frac{z}{t}}t^{-c}.\eeq
\qed

\noindent
{\bf ${}_0F_1$ function} 
\[F(c;z):=\sum_{j=0}^\infty
\frac{1}{
(c)_j}\frac{z^j}{j!}.\]
{\bf Olver normalized  ${}_0F_1$ function}:
\[ {\bf F}  (c;z):=\frac{F(c;z)}{\Gamma(c)}=
\sum_{j=0}^\infty
\frac{1}{
\Gamma(c+j)}\frac{z^j}{j!}.\]

\noindent{\bf Basic pair of recurrence relations:}
\begin{align}
  \p_z\mathbf{F}(c;z)&=\mathbf{F}(c+1;z),\\
  ( \p_z+c)\mathbf{F}(c+1;z)&=\mathbf{F}(a;c;z).
                                             \end{align}

                                             \noindent
                                             {\bf Integral representation} for all parameters:
\begin{eqnarray*}
\frac{1}{2\pi \i}\int\limits_{]-\infty,0^+,-\infty[}
\e^t\e^{\frac{z}{t}}t^{-c}\d t
&=& {\bf F}  (c;z),\ \ \ \ \Re z>0.\end{eqnarray*}

\noindent
{\bf Chebyshev solutions} for $k=0,1,2,\dots$:
\begin{align}
{\bf  F}(\tfrac32+k;z)&=\frac{2^{2k}}{\sqrt\pi}\partial_z^k\frac{\sinh 
                   2\sqrt{z}}{\sqrt{z}},\\
  {\bf  F}(\tfrac12-k;z)&=\frac{z^{\frac12+k}}{\sqrt\pi}\partial_z^k\frac{\cosh 
                   2\sqrt{z}}{\sqrt{z}}.
  \end{align}

\subsection{Hermite equation}

{\bf Hermite operator:}
\beq\S({a}):=\p_z^2-2z\p_z-2{a}.\eeq
\noindent {\bf Parameters:} 
\begin{align*}
  \sigma(z)&=1,\quad
  \kappa(z)=-2z,\quad
  \omega=-2a+1.\end{align*}
{\bf Weight:}
\[  \rho(z)=\e^{-z^2}.\]
\text{\bf Basic symmetry:}\[\e^{-z^2}\cS(a;z,\p_z) \e^{z^2}=
-\cS(1-a;w,\p_w),\qquad
w=\pm \i z.\]
\bet[about integral representations]\ben\item
Let $[0,1]\ni t\mapsto\gamma(t)$ satisfy
\[\e^{t^2}(t-z)^{-{a}-1}\Big|_{\gamma(0)}^{\gamma(1)}=0.\]
Then
\beq\S({a})\quad\text{annihilates}\quad
\int_\gamma\e^{t^2}(t-z)^{-{a}}\d t.\label{dad10}\eeq
\item
Let $[0,1]\ni t\mapsto\gamma(t)$ satisfy
\[\e^{-t^2-2zt}t^{{a}}\Big|_{\gamma(0)}^{\gamma(1)}=0.\]
Then\beq
\S({a})
\quad\text{annihilates}\quad\int_\gamma
\e^{-t^2-2zt}t^{a-1}\d t.\label{dad11}\eeq
\een\label{dad12}\eet

\proof
We check that 
\begin{eqnarray*}
\S({a})\,\e^{t^2}(t-z)^{-{a}}=&-{a} \p_t
\e^{t^2}(t-z)^{-{a}-1},\\
\S({a})\,
\e^{-t^2-2zt}t^{a-1}=&-2\p_t\e^{-t^2-2zt}t^{{a}}.
\end{eqnarray*}
\qed\\
{\bf Hermite function:} 
\begin{eqnarray}\label{hermite}
S(a;z)&:=&
z^{-{a}}F\Big(\frac{{a}}{2},\frac{{a}+1}{2};-;-z^{-2}\Big).
\end{eqnarray}

\noindent{\bf Basic pair of recurrence relations:}
\begin{align}
  \p_zS(a;z)&=-aS(a+1;z),\\
  \big( \p_z-2z\big)S(a+1;z)&=S(a;z).
                                             \end{align}

\noindent
{\bf Integral representation} for $z\not\in]-\infty,0]$.
  For $\red 0<{\Re}a$:
\begin{eqnarray*}
\int\limits_0^\infty\e^{-t^2-2tz}t^{a-1}\d t
&=&2^{-a}\Gamma(a)S(a;z);
\end{eqnarray*} and for all parameters:
\begin{eqnarray*}
-\i\int\limits_{]-\i\infty,z^-,\i\infty[}\e^{t^2}
(z-t)^{-a}
\d
  t 
&=&\sqrt\pi S(a;z).
\end{eqnarray*}

\subsection{Lie algebras of symmetries type by type}

Miller's Lie algebra is responsible only for one pair of recurrence
relations of hypergeometric class functions. {\red Some of  types from the
hypergeometric class possess larger sets of recurrence
relations, which were
described e.g. in \cite{De1,De2}. These recurrence relations can be
interpreted in terms of roots of  larger Lie algebras
of symmetries, described e.g. in \cite{De2}.
}

Miller's Lie algebra possesses a Weyl group isomorphic to $\zz_2$.
Some types of hypergeometric class equations possess larger groups of discrete symmetries.

In the following table we list all nontrivial types of hypergeometric
class equations. {\red We include
the {\em Gegenbauer equation}
\beq \Big((1-w^2)\partial_w^2-2(1+\alpha )w\partial_w
+\lambda ^2-\big(\alpha +\tfrac{1}{2}\big)^2\Big)f(w)=0,\label{gege0}
\eeq
treating it as a seperate type. Strictly speaking, the Gegenbauer type
is contained in the
  ${}_2F_1$ type. In fact, the transformation $w\mapsto
  \frac{1-w}{2}$ transforms the Gegenbauer equation into a special case
  of the ${}_2F_1$ equation.
The Gegenbauer equation has the special property of} the {\em mirror
symmetry} $w\mapsto-w$. In the remaining part of this paper
  there is  no need to consider it separately, however it has special
  symmetry properties.

  The second and third column are based on Section \ref{Miller's Lie
    algebra} of this paper. {\red In particular, the second column can be
  compared with \eqref{miller1}--\eqref{miller3}.}

The fourth and fifth column are based on
  \cite{De2}. They describe a more complete Lie algebra of symmetries and
  the corresponding group of discrete symmetries.

\medskip

\begin{tabular}{lllll}
  Equation&Miller's Lie algebra&surface $\Omega$&Lie algebra \cite{De2}& Discrete
                                                              symmetries
                                                                \cite{De2}\\[1ex]
  \hline\\
  ${}_2F_1$&$sl(2,\cc)\oplus\cc$&sphere&$so(6,\cc)$&symmetries of cube\\[2ex]
             ${}_1F_1$&$osc(\cc)$&paraboloid&$sch(2,\cc)$&$\zz_2\times\zz_2$\\[2ex]
                        ${}_2F_0$&$sl(2,\cc)\oplus\cc$&null quadric&$sch(2,\cc)$&$\zz_2\times\zz_2$\\[2ex]
                                   ${}_0F_1$&$\cc^2\rtimes so(2,\cc)\oplus\cc$&paraboloid&$\cc^2\rtimes so(2,\cc)$&$\zz_2$\\[2ex]
                                              Gegenbauer&$sl(2,\cc)\oplus\cc$&sphere&$so(5,\cc)$&
  symmetries of                                                                                         square\\[2ex]
                                                          Hermite&$osc(\cc)$&cylinder&$sch(1,\cc)$&$\zz_4$.
  \end{tabular}
  \medskip

 Above,  $sch(n,\cc)$
  denotes the complex Schr\"odinger Lie algebra, the Lie algebra of
  generalized symmetries of the heat equation in dimension $n$---see
  e.g. \cite{De2}.

{\red The group from the fifth column  always  contains various
symmetries described in our paper. First of all, all types possess
the basic symmetry, which is the generator of 
the Weyl group of Miller's Lie algebra. Some of the types possess
additional symmetries.

The power symmetry is a property of the ${}_2F_1$, ${}_1F_1$ and $F_1$
equation. For the $F_1$ equation it coincides with the basic symmetry.

The inversion symmetry is a feature of the ${}_2F_1$ equation, and
connects  the ${}_1F_1$ equation with the  ${}_2F_0$ equation.

The ${}_2F_1$ and ${}_2F_0$ equations are invariant wrt swapping $a$ and
$b$.

The Gegenbauer and Hermite equation are invariant wrt the
symmetry $z\mapsto-z$.

The Lie algebra of the fourth 
column always contains the corresponding Miller's Lie algebra (possibly, without the 
trivial term 
$\oplus\cc$). By applying discrete symmetries from the fifth column to
Miller's Lie algebra, we can enlarge Miller's Lie algebra to a
subalgebra of  the  Lie algebra from the
fourth column.

}

\init
\section{Hypergeometric  polynomials}

\label{Hypergeometric  polynomials}

Polynomial solutions of hypergeometric class equations will be called {\em  
  hypergeometric class polynomials}. 
In this 
section we describe these solutions in detail.

We have already mentioned in \eqref{ladder2a} that there exist 
ladders of solutions of hypergeometric class equations consisting of 
polynomials. 
Note that according to our conventions, these ladders are {\em descending}.
This is related to the usual convention for parameters of various
types of hypergeometric functions: in order to get a 
polynomial of degree $n\in\nn_0$, the parameter 
$a$  takes the value $-n$.

In Subsections \ref{Differential operators with polynomial
  eigenfunctions}
and
\ref{s.12} we describe the algebraic theory of hypergeometric class 
polynomials. They are centered around the so-called (generalized)
Rodrigues formula.

Under some conditions on $\sigma$ and $\kappa$
hypergeometric class 
polynomials can be viewed 
as eigenfunctions of a certain self-adjoint Sturm-Liouville
operator acting on an appropriate weighted Hilbert space.
Besides, they form an  orthogonal basis of this Hilbert space.
To explain this point of view, we devote
Subsections \ref{Orthogonal polynomials},
\ref{Hermiticity of Sturm-Liouville operators}
and \ref{Selecting endpoints for  Sturm-Liouville operators}
to a few general remarks 
about orthogonal polynomials and Sturm-Liouville
operators. In particular, we find some useful  conditions for the
Hermiticity of Sturm-Liouville operators. These conditions
yield  the weight function $\rho$, and also specify
possible endpoints of the interval $]a,b[$.

{\red In the remaining subsections
 we return to the theory of hypergeometric class
polynomials. We are convinced that the unified approach to their
theory, which we present,
based partly on \cite{NU},  has a considerable pedagogical value.}

\subsection{2nd order differential operators with polynomial eigenfunctions}
\label{Differential operators with polynomial eigenfunctions}

The following well-known and easy proposition shows  that 
hypergeometric class operators have many polynomial eigenfunctions. Actually, it is the property that characterizes this class among all
second order
differential operators.

\bep  Let 
$\sigma(z)$, $\tau(z)$, $\mu(z)$ be arbitrary functions.
Let $P_n(z)$,
$n=0,1,2$, be
polynomials such that $\deg P_n=n$ and
 $\eta_n\in\cc$. Suppose that
\[\left(\sigma(z)\p_z^2+\tau(z)\p_z+\mu(z)+\eta_n\right)P_n(z)=0,\ \ n=0,1,2.\]
Then $\sigma(z)$ is a polynomial of degree at most $2$, $\tau(z)$
 is a
 polynomial of degree at most 1 and $\mu(z)$ a polynomial of degree
 $0$ (a complex number).
\label{con11}\eep

Note that for differential operators from Prop. \ref{con11}, by
replacing $\mu+\eta_n$ with $\eta_n$, without limiting the generality
we can assume that $\mu=0$.

 Let us describe some other simple facts about polynomial solutions of
hypergeometric class equations:

{\red
\bep Let $\deg\sigma\leq2$,  $\deg\tau\leq1$, $\eta\in\cc$,
\begin{enumerate}
\item Suppose that $P$ is a polynomial of degree $n$  solving
\beq\label{propio}
\left(\sigma(z)\p_z^2+\tau(z)\p_z+\eta\right)P(z)=0.\eeq
Then
\beq\begin{array}{l}
n(n-1)\frac{\sigma''}{2}+n\tau'+\eta=0.\end{array}
\label{con1}\eeq 
\item
If
\beq\label{contra} k\frac{\sigma''}{2}+\tau'\neq0,\quad k\in\nn_0,\eeq
then the space of polynomial solutions of \eqref{propio}
is at most $1$-dimensional.\label{propio0}\end{enumerate}
\eep

\proof
 Differentiating $n$ times \eqref{propio} we obtain
\eqref{con1}$\times P^{(n)}=0$. This implies 1.}

Suppose the space of polynomial
solutions of  \eqref{propio} is 2-dimensional. We can assume that the
degrees of these solutions are $0\leq n_1<n_2$. By \eqref{con1},
\beq\label{popo}
n_i(n_i-1)\frac{\sigma''}{2}+n_i\tau'+\eta=0,\qquad i=1,2.\eeq
Subtracting \eqref{popo} for 2 and 1 and dividing by $n_1-n_2$ we obtain
\beq
(n_1+n_2-1)\frac{\sigma''}{2}+\tau'=0.\eeq
Now possible values of $k:=n_1+n_2-1$ are $0,1,2,\dots$. This
contradicts \eqref{contra}.
\qed

\subsection{Raising and lowering operators. Rodriguez formula}
\label{s.12}

The following theorem shows how to construct   hypergeometric class polynomials
and describes some of their properties, which are easy to describe in a
unified fashion. As usual, $\sigma$ denotes a polynomial of degree
$\leq2$, $\kappa$ a polynomial of degree $\leq1$.

\bet
\begin{enumerate}\item Suppose that $\omega\in\cc$ and 
a polynomial of degree $n$ solves
\begin{align}
\label{pso}  0&=\Big({\red \cH}(\sigma,\kappa)+\omega \Big)P(z)
\\
&=\Big(\sigma(z)\p_z^2+(\sigma'(z)+\kappa(z))\p_z +\frac{\kappa'}{2}+\omega\Big)P(z).
\end{align}
Then
\beq\label{eigeno}
n(n+1)\frac{\sigma''}{2}+(n+\tfrac12)\kappa'+\omega=0.\eeq
\item Suppose that
  \beq k\frac{\sigma''}{2}+\kappa'\neq0,\quad k=2,3,\dots.\eeq
Then for any $\omega\in\cc$ the space of polynomial  solutions of \eqref{pso} is at
most $1$-dimensional.
\item
 Define
\begin{align}\label{rod}
P_n(\sigma,\kappa;z)&:=\frac{1}{n!}\rho^{-1}(z)\p_z^n\sigma^n(z)\rho(z)\\
&=\frac{1}{2\pi \i}\rho^{-1}(z)\int\limits_{[z^+]}
\sigma^n(z+t)\rho(z+t)t^{-n-1}\d t, \label{rod2}
\end{align}
where $[z^+]$ denotes a loop around $z$ in the counterclockwise direction.
Then $P_n(\sigma,\kappa;z)$ is a polynomial of degree $n$ or less, 
and we have
\begin{align}                       \label{pop}
0&=\Big({\red \cH}(\sigma,\kappa)-n(n+1)\frac{\sigma''}{2}-(n+\12)\kappa'\Big)
P_n(\sigma,\kappa;z)\\
&=\Big(\sigma(z)\p_z^2+(\sigma'(z)+\kappa(z))\p_z
-n(n+1)\frac{\sigma''}{2}-n\kappa'\Big)P_n(\sigma,\kappa;z).
\notag
\end{align}
\end{enumerate}\eet

\proof 1. and 2. 
are just a reformulations of Proposition  \ref{propio0}, where 
instead of parameters $\tau,\eta$ we use $\kappa,\omega$. 

To show 3. we note that
the integral representation \eqref{rod2} satisfies the assumptions
of Theorem \ref{in1} about Euler 
transforms. This representation  implies both \eqref{pop} and
\eqref{rod}.
\qed

\eqref{rod}
is  usually called 
the {\em Rodriguez formula}.

Consider the descending ladder starting at
$\kappa_0(z):=\kappa(z)$, $\omega_0=-\frac{\kappa'}{2}$, that is,
\begin{align}
  \kappa_{-n}(z)&:=\kappa(z)-n\sigma'(z),\\
  \omega_{-n}&:=-(n+\tfrac12)\kappa'+n^2\frac{\sigma''}{2}=-\kappa_{-n}'(n+\tfrac12)-n(n+1)\frac{\sigma''}{2}.
\end{align}

\bep The polynomials $P_n(\sigma,\kappa_{-n};z)$ satisfy the corresponding
hypergeometric class equation
\begin{align}
\Big({\red \cH}(\sigma,\kappa_{-n})+\omega_{-n}\Big) 
  P_n(\sigma,\kappa_{-n};z)&=0,\end{align}
the recurrence relations\begin{align}
(\sigma(z)\p_z+\kappa_{-n}(z))P_n(\sigma,\kappa_{-n};z)&=
(n+1)P_{n+1}(\sigma,\kappa_{-n-1};z),\label{recur1}
\\
\p_zP_{n+1}(\sigma,\kappa_{-n-1};z)&=
\Big(n\frac{\sigma''}{2}+\kappa_{-n}')p_{n}(\sigma,\kappa_{-n};z);\label{recur2}\end{align}
and have a generating function
\begin{align}\label{recur3}
  \frac{\rho(z+t\sigma(z))}{\rho(z)}&=\sum\limits_{n=0}^\infty 
t^nP_n(\sigma,\kappa_{-n};z).\end{align}
\eep

\subsection{Orthogonal polynomials}
\label{Orthogonal polynomials}

Consider
$-\infty\leq a<b\leq+\infty$.
Suppose that $]a,b[\ni x\mapsto \rho(x)$ is a positive measurable function. Define the Hilbert space
\begin{align}
  L^2(]a,b[,\rho)&:=\Big\{f\quad \text{measurable on}\quad ]a,b[ \ |\
  \int_a^b\rho(x)|f(x)|^2\d x<\infty\Big\},\\
\text{with the scalar product }\quad
                 &(f|g):=\int_a^b\overline{f(x)}g(x)\rho(x)\d
                   x.\end{align}

Assume in addition
\beq\int_a^b\rho(x)|x|^n\d x<\infty,\quad n\in\nn.\eeq
Then the space of polynomials is contained in $ L^2(]a,b[,\rho)$.
Applying the Gram-Schmidt orthogonalization to the sequence
$1,x,x^2,\dots$ we can define an orthogonal family of polynomials
$p_0,p_1,p_2,\dots$ The following simple criterion is proven e.g. in
\cite{ReedSimon}:
\bet\label{polypo}
Suppose that for some $\epsilon>0$
\beq\int_a^b\e^{\epsilon|x|}\rho(x)\d x<\infty.\eeq
Then polynomials are dense in $ L^2(]a,b[,\rho)$. Therefore,
the family $p_0,p_1,p_2,\dots$ is an orthogonal basis of
$ L^2(]a,b[,\rho)$. \eet

\subsection{Hermiticity of Sturm-Liouville operators}
\label{Hermiticity of Sturm-Liouville operators}

In most of this paper we avoid using functional analysis. However, for
hypergeometric class polynomials we will make an exception. In fact,
often it is natural to view them as  orthogonal bases consisting of
eigenfunctions of certain self-adjoint 
operators.

Let us
briefly recall some elements of the theory of operators on Hilbert spaces.
  Let $\cH$ be a Hilbert space with the scalar product
  $(\cdot|\cdot)$. Let
  $\cA$  be an operator on the domain $\cD\subset\cH$. We say that $\cA$
  is {\em Hermitian } on $\cD$ in the sense of $\cH$ 
if 
\beq (f|\cA g)=(\cA f|g),\quad f,g\in\cD.\eeq
(Unfortunately, in most  of the contemporary  mathematics 
literature, instead of {\em Hermitian} the word  {\em symmetric} is used,
which is a confusing misnomer).
We say that $\cA$ is {\em self-adjoint} if it is Hermitian and its spectrum
 is real. We say that it is {\em essentially self-adjoint} if it has a
 unique self-adjoint extension.

Let $\cA$ be Hermitian. Suppose that $f_i\in \cD$, $i=1,2$, and
 \beq \cA f_i=\eta_i f_i,\quad f_i\in \cD,\quad i=1,2.\eeq
 Then it is easy to see that $\eta_i$ are real, and if
 $\eta_1\neq\eta_2$, then $(f_i|f_j)=0$. Therefore, eigenfunctions of
 Hermitian operators can be arranged into orthogonal families.
If an operator is Hermitian and possesses an orthogonal
 basis of eigenvectors, then it is  is essentially
 self-adjoint on (finite) linear combinations of its eigenvectors.

Second order differential operators on a segment of the real line are
often called {\em Sturm-Liouville operators}. 
Let us make some remarks about the general theory of such operators,
which will be useful in the analysis of orthogonality properties of
hypegeometric class polynomials.

Consider
$-\infty\leq a<b\leq+\infty$. Suppose that
$\sigma,\rho$ are  functions on
$]a,b[$ (not necessarily polynomials).
Consider an operator
of the form
\beq \cA:=\sigma(x)\partial_x^2+\tau(x)\partial_x\label{herm1}\eeq 
(Thus we consider an arbitrary Sturm-Liouville operator without the 
zeroth order term).

Let $\rho$ satisfy 
\beq\label{rro}
\sigma(x)\rho'(x)=(\tau(x)-\sigma'(x))\rho(x).\eeq 
Often it is convenient to rewrite \eqref{herm1} in the  form 
\beq\label{sturm}
\cA=\rho(x)^{-1}\partial_x\rho(x)\sigma(x)\partial_x.\eeq

Suppose that $\sigma,\tau$ are real and $\rho$  is positive.
It is useful to interpret $\cA$ as an (unbounded) operator on the
Hilbert space
$  L^2(]a,b[,\rho)$.
Indeed, the following easy theorem shows that $\cA$ 
is Hermitian on smooth functions that vanish close to the endpoints of
the interval.

\bet Let
\beq
\cD_0=\{f\in C^\infty(]a,b[)
\ :\ f=0\ \ \hbox{in a neighborhood of}\ \ a,b\}.\label{endpo}
\eeq
Then $\cA$ is Hermitian  on $\cD_0$ in the sense of the Hilbert space
$L^2(]a,b[,\rho)$.
\label{herm3}\eet

\subsection{Selecting endpoints for  Sturm-Liouville operators}
\label{Selecting endpoints for  Sturm-Liouville operators}

Unfortunately, the operator $\cA$  is rarely essentially self-adjoint on \eqref{endpo}.
We will see that  $\cA$ 
 is still Hermitian on functions that do not vanish near the
endpoints,  if these endpoints satisfy appropriate conditions:

\bet 
Let $-\infty<a<b<+\infty$ and
\[\sigma(a)\rho(a)=\sigma(b)\rho(b)=0.\]
Then $\cA$ is Hermitian on the domain $C^2([a,b])$ in the sense of the space
$ L^2(]a,b[,\rho)$.
\label{herm4}\eet

\proof Let $g,f\in C^2([a,b])$.
\begin{eqnarray*}
(g|\cA f)&=&
\int_a^b\rho(x)\bar g(x)\rho(x)^{-1}
\partial_x\sigma(x)\rho(x)\partial_x f(x)\d x\\
&=&
\int_a^b\bar g(x)
\partial_x\sigma(x)\rho(x)\partial_x f(x)\d x\\
&=&
\bar{g(x)}\rho(x)\sigma(x) f'(x)\Big|_a^b-\int_a^b(\partial_x  \bar g(x))
\sigma(x)\rho(x)\partial_x f(x)\d x\\
&=&
-\bar{g'(x)}\rho(x)\sigma(x)f(x)\Big|_a^b+\int_a^b(\partial_x
  \rho(x)\sigma(x)\partial_x\bar{g(x)})f(x)\d x\\
&=&
\int_a^b\rho(x)\bar{\left(\rho(x)^{-1}\partial_x
\sigma(x)\rho(x)\partial_x
  g(x)\right)} 
f(x)\d x=(\cA g|f).
\end{eqnarray*}
\qed

Analogously we prove the following fact:

\bet 
Let
\[\lim_{x\to-\infty}\sigma(x)\rho(x)|x|^n=
\lim_{x\to+\infty}\sigma(x)\rho(x)|x|^n=
0,\ \ \ n\in\nn.\]
Then $\cA$ is Hermitian on the domain consisting of polynomial
functions in the sense of the Hilbert space 
$ L^2(]-\infty,\infty[,\rho)$.
\label{herm4a}\eet

Obviously,  statements similar to Theorems \ref{herm4} and 
\ref{herm4a}  hold if $a=-\infty$ and $b$ is finite, or $a$ is finite
and $b=\infty$.

We will see later that
in applications to hypergeometric class equations the conditions of Theorems
\ref{herm4} and 
\ref{herm4a}  are often sufficient to guarantee the essential
self-adjointness
of the  operator $\cA$ on the set of polynomials.  

\subsection{Orthogonality of hypergeometric  class polynomials}
\label{s.13}

In some cases hypergeometric class polynomials can be interpreted as an
orthogonal basis of a certain weighted Hilbert space. They are then
often called {\em (very) classical orthogonal polynomials}.

Let us fix $\sigma$, $\kappa$ and $\rho$ satisfying (\ref{e1b}).
{\red Let $\sigma,\kappa$ be real and $\rho$ positive.}
Choose
 an interval $]a,b[$, where
$-\infty\leq a<b\leq+\infty$ satisfy the conditions of
Theorem \ref{herm4} or \ref{herm4a}.
By these theorems, the operator
\begin{align}\red \cH(\sigma,\kappa)
  &:=\sigma(x)\partial_x^2+\big(\kappa(x)+\sigma'(x)\big)\partial_x+\frac{\kappa'}{2}=\rho(x)^{-1}\partial_x\rho(x)\sigma(x)\partial_x+\frac{\kappa'}{2}.\end{align}
is Hermitian on the space of polynomials in the sense of the Hilbert space
$ L^2(]a,b[,\rho)$. Therefore, 
polynomial eigenfunctions of $\cH(\sigma,\kappa)$ are
pairwise orthogonal (at least if the corresponding eigenvalues are
distinct, which is generically the case, see \eqref{eigeno}).
The following theorem says more: we compute the square norm of the
polynomials obtained from the Rodriguez formula.

\bet
\begin{align}
  \int_a^bP_m(\sigma,\kappa;x)P_n(\sigma,\kappa;x)\rho(x)\d x
=\frac{\delta_{mn}}{n!}&\prod_{j=1}^n\Big(-\kappa'+j\frac{\sigma''}{2}\Big)
\int_a^b\sigma^n(x)\rho(x)\d x.\label{pq1}
\end{align}
\label{orth}\eet

\proof
It is enough to assume that $m\leq n$. Write $P_n(x)$ for $P_n(\sigma,\kappa;x)$.
\begin{align}
  &\int_a^bP_m(x)P_n(x)\rho(x)\d x\\
  =&\frac{1}{n!}\int_a^bP_m(x)\partial_x^n\sigma^n(x)\rho(x)\d x\\
    =&\frac{(-1)^n}{n!}\int_a^b\Big(\partial_x^nP_m(x)\Big)\sigma^n(x)\rho(x)\d x.
\end{align}
This is zero if $m<n$ and for $m=n$ we use
\beq
\partial_x^nP_n(x)=\prod_{j=1}^n\Big(\kappa'-j\frac{\sigma''}{2}\Big).\eeq
\qed

\subsection{Review of types of hypergeometric class polynomials}

In the remaining part of the section we review various types of
hypergeometric class polynomials. We will discuss only the properties that
follow directly from the general theory described in  previous
subsections.
Here is the list of items that we will cover:
\begin{itemize}
\item The choice of $\sigma,\kappa$ and the corresponding weight
  $\rho$.
  \item The Rodriguez formula \eqref{rod} defining the polynomial and
    its expression in
    terms of hypergeometric class functions.
  \item The degree of the polynomial, see \eqref{con1}.
  \item The differential equation, see \eqref{pop}.
    \item The pair of recurrence relations that follow from the
      Rodriguez formula, see \eqref{recur1} and \eqref{recur2}.
      \item The generating function related to the Rodriguez formula,
        see \eqref{recur3}.
        \item The choice of endpoints for which the  corresponding
          Sturm-Liouville operator is essentially self-adjoint on the
          space of polynomials (if applicable), see Theorems
 \ref{herm4} ,
\ref{herm4a} and \ref{polypo}.       \item The square norm (if applicable), see
Theorem \ref{orth}.
          \end{itemize}

          Hypergeometric class polynomials
          possess various useful features that will not be listed in
          the following subsections. For instance, they typically
          have additional recurrence relations and generating
          functions.
          We list only those that follow directly from the general
          theory described above.

\subsection{Jacobi polynomials}

Consider  $\alpha,\beta\in\cc$ and
\beq
\sigma(z)=1-z^2,\quad\kappa(z)=\alpha(1-z)+\beta(1+z),\quad
\rho(z)=(1-z)^\alpha(1+z)^\beta.\eeq
For $n\in\{0,1,\dots\}$ set
\begin{eqnarray}\label{jaco}
P_n^{\alpha,\beta}(x)&=&
\frac{(-1)^n}{2^nn!}
(1-x)^{-\alpha}(1+x)^{-\beta}\p_x^n
(1-x)^{\alpha+n}(1+x)^{\beta+n}\\\label{holso1}
&=&\frac{(1+\alpha)_n}{n!}
{}_2F_1\Big(-n,n+\alpha+\beta+1;\alpha+1;\frac{1-x}{2}\Big)
.\end{eqnarray}
Then
$P_n^{\alpha,\beta}$ is a polynomial of degree at most $n$. More precisely:
\ben\item If $\alpha+\beta\not\in\{-2n,\dots,-n-1\}$, then
$\deg P_n^{\alpha,\beta}= n$. It is then up
to a coefficient the  unique eigenfunction of the operator
\beq \red \cH(\sigma,\kappa)=(1-x^2)\p_x^2+(\beta-\alpha-
(\alpha+\beta+2)x)\p_x+\frac{\beta-\alpha}{2}\eeq
which is a polynomial of degree $n$.
\item If  $\alpha+\beta\in\{-2n,\dots,-n-1\}$,  but 
  $\alpha\not\in\{-n,\dots,-1\}$ (or, equivalently, 
    $\beta\not\in\{-n,\dots,-1\}$), then 
$\deg P_n^{\alpha,\beta}= -\alpha-\beta-n-1$.
\item If  $\alpha+\beta\in\{-2n,\dots,-n-1\}$,  but 
  $\alpha\in\{-n,\dots,-1\}$ (or, equivalently, 
    $\beta\in\{-n,\dots,-1\}$), then 
$P_n^{\alpha,\beta}= 0$.\een

$P_n^{\alpha,\beta}$ satisfy the  Jacobi equation, which is a
slightly modified ${}_2F_1$  equation
\[\left(
(1-x^2)\p_x^2+(\beta-\alpha-
(\alpha+\beta+2)x)\p_x+n(n+\alpha+\beta+1)\right)P_n^{\alpha,\beta}(x)=0.
\]
Recurrence relations:
\begin{eqnarray}
\partial_xP_n^{\alpha,\beta}(x)&=&\frac{\alpha+\beta+n+1}{2}
P_{n-1}^{\alpha+1,\beta+1} ,\label{jaco2}\\
-\frac{(1-x^2)\partial_x+\beta-\alpha-(\alpha+\beta)x
}{2}
P_n^{\alpha,\beta}(x)&=&(n+1)P_{n+1}^{\alpha-1,\beta-1}(x).\label{rod4b}
\end{eqnarray}
Generating function:
\beq
\sum_{n=0}^\infty
P_n^{\alpha-n,\beta-n}(x)2^nt^n
=\bigl(1+t(1+x)\bigr)^\alpha\bigl(1-t(1-x)\bigr)^\beta.\label{wiel4}\eeq

If $\alpha,\beta>-1$, then $\cH(\sigma,\kappa)$ is self-adjoint on the space of
polynomials  in the sense of
\[L^2(]-1,1[,(1-x)^\alpha(1+x)^\beta).\]
Jacobi polynomials are its eigenfunctions and form an orthogonal
basis with the normalization
\beq\int_{-1}^1
P_n^{\alpha,\beta}(x)^2(1-x)^\alpha(1+x)^\beta\d x
=\frac{\Gamma(1+\alpha+n)\Gamma(1+\beta+n)2^{\alpha+\beta+1}}
{(1+2n+\alpha+\beta)n!\Gamma(1+\alpha+\beta+n)}
.\label{ppw1}
\eeq

\subsection{Laguerre polynomials}

Consider   $\alpha\in\cc$ and
\beq
\sigma(z)=z,\quad\kappa(z)=\alpha-z,\quad\rho(z)=\e^{-z}z^\alpha.\eeq
 For $n\in\nn$ set
\begin{eqnarray*}
L_n^\alpha(x)&=&\frac{1}{n!}\e^xx^{-\alpha}\p_x^n\e^{-x}x^{n+\alpha}\\
&=&\frac{(1+\alpha)_n}{n!}{}_1F_1(-n;1+\alpha;x).
\end{eqnarray*}
Then $L_n^\alpha$ is a polynomial of degree $n$.
It is a unique (up to a coefficient) eigenfunction of the operator
\beq\red
\cH(\sigma,\kappa)=x\p_x^2+(\alpha+1-x)\p_x-\frac12\eeq
which is a polynomial of degree $n$. 
$L_n^\alpha$ satisfy the Laguerre equation, which is the ${}_1F_1$
equation with modified parameters:
\[\left(
x\p_x^2+(\alpha+1-x)\p_x+n\right)L_n^\alpha(x)=0.\]
Recurrence relations:
\begin{eqnarray}
(x\partial_x+\alpha-x)L_n^\alpha(x)&=&(n+1)L_{n+1}^{\alpha-1}(x),\\
\partial_x L_n^\alpha(x)&=&-L_{n-1}^{\alpha+1}(x).\label{rod4a}
\end{eqnarray}
Generating function:
\beq
\e^{-tz}(1+t)^{\alpha}
=\sum\limits_{n=0}^\infty t^n L_n^{\alpha-n}(z).
\eeq
If $\alpha>-1$, then $\cH(\sigma,\kappa)$ is essentially self-adjoint on the space of
polynomials in the sense of $L^2([0,\infty[,\e^{-x}x^\alpha)$.
 Laguerre polynomials are its eigenfunctions and form an orthonormal
basis 
with the normalization
\[\int_0^\infty L_n^\alpha
(x)^2 x^\alpha\e^{-x}\d x=\frac{\Gamma(1+\alpha+n)}{n!} .\]

\subsection{Bessel polynomials}
\label{Bessel polynomials}
Consider $\theta\in\cc$ and
\beq
\sigma(z)=z^2,\quad\kappa(z)=-1+\theta 
                           z,\quad\rho(z)=\e^{-z^{-1}}z^\theta.\eeq
For $n=0,1,\dots$ set
\begin{align}
B_n^\theta(z):&=\frac{1}{n!}
                z^{-\theta}\e^{z^{-1}}\p_z^n\e^{-z^{-1}}z^{\theta+2n}\\
              &=\frac{1}{n!}{}_2F_0(-n,n+\theta+1;-;z)\\
  &=(-z)^nL_n^{-\theta-2n-1}(-z^{-1}).\end{align}
Then $B_n^\theta$ is a polynomial of degree $n$.
It is a unique (up to a coefficient) eigenfunction of the operator
\beq
\red\cH(\sigma,\kappa)=z^2\p_z^2+(-1+(2+\theta)z)\p_z+\frac\theta2\eeq
which is a polynomial of degree $n$. 
$B_n^\theta$ satisfy the ${}_2F_0$ equation with adjusted parametes:
\begin{eqnarray*}
\Big(z^2\p_z^2+(-1+(2+\theta)z)\p_z-\frac12n(1+\theta+n)\Big)B_n^\theta(z)&=&0.
\end{eqnarray*}
Recurrence relations:
\begin{eqnarray*}
\p_zB_n^\theta(z)&=&-(n+\theta+1)B_{n-1}^{\theta+2}(z),\\
\left(z^2\p_z-1-\theta z\right)B_n^\theta(z)&=&-(n+1)B_{n+1}^{\theta-2}(z).
\end{eqnarray*}
Generating function: 
\beq 
(1+tz)^{\theta}\exp\Big(\frac{-t}{1+tz}\Big)=\sum\limits_{n=0}^\infty 
t^nB_n^{\theta-2n}(z).\eeq

Bessel polynomials do not form an orthogonal basis on any interval.

\subsection{Hermite polynomials}

Consider
\beq
\sigma(z)=1,\quad\kappa(z)=-2z,\quad\rho(z)=\e^{-z^2}.\eeq
For $n=0,1,\dots$ set
\begin{align}
  H_n(x)&=\frac{(-1)^n}{n!}\e^{x^2}\p_x^n\e^{-x^2}\\
        &=\frac{2^n}{n!}S(-n;x),\end{align}
      where $S(a,x)$ is the Hermite function defined in \eqref{hermite}.
Then $H_n$ is a polynomial of degree $n$ 
and is (up to a multiplicative constant) the only eigenfunction  
of the operator \beq
\red\cH(\sigma,\kappa)=\p_x^2-2x\p_x-1\eeq which is a polynomial of degree $n$. 
It satisfies the  Hermite equation
\[(\p_x^2-2x\p_x+2n)H_n(x)=0.\]
Recurrence relations:
\begin{align}
(-\partial_x+2x)H_n(x)&=(n+1)H_{n+1}(x),\\
  \partial_x H_n(x)&=2H_{n-1}(x).\label{rod4}\end{align}
Generating function:
\begin{align}\sum_{n=0}^\infty t^n H_n(x)&=\e^{2tx-t^2}.\label{wiel6}
\end{align}

The operator $\cH(\sigma,\kappa)$ is essentially self-adjoint on the space of
polynomials in the sense of
$L^2(\rr,\e^{-x^2})$. Hermite polynomials are its eigenfunctions
and form an orthogonal basis with the normalization
\[\int_{-\infty}^\infty H_n(x)^2\e^{-x^2}\d x=\frac{\sqrt\pi2^n}{n!}.\]

\begin{remark} The definition of Hermite polynomials that we gave is
consistent with the generalized Rodrigues formula \eqref{rod}.
In the literature one can also find other conventions for Hermite
polynomials, e.g.
$H_n(x):=(-1)^n\e^{x^2}\p_x^n\e^{-x^2}.$
\end{remark}

\init
\section{Degenerate case}
\label{Degenerate case}

If two indices of  a regular-singular (also called Fuchsian) point of a
differential equation differ by an integer, then the usual Frobenius
method  \cite{WW}
produces  in general, up to a coefficient, only one  solution. We call this
case {\em degenerate}. In this section we will be discuss the
degenerate case of hypergeometric class equations. Without limiting
the generality, the regular-singular point will be $0$.

\subsection{Unified hypergeometric function in the degenerate case}

Suppose that we are in the setting of Section \ref{sec.sec}. That means, $\sigma(0)=0$,
 $\sigma'(0)=1$, and we set $\kappa(0)=m$. Thus, 
\beq\label{param}
\sigma(z)=\frac{\sigma''}{2}z^2+z,\quad\kappa(z)=m+\kappa' z,\quad
\omega\eeq
are the parameters of the equation we consider, see
\eqref{olver.1}. As in \eqref{olver.2}, we also introduce the
transformed parameters:
\beq
\kappa^\circ(z):=-m+(\kappa'-m\sigma'')z,\quad\omega^\circ:=\omega-m\kappa'+m^2\frac{\sigma''}{2},\quad
m^\circ=-m.\eeq
We assume in addition that $m$ is an integer.
Using
\beq\frac{1}{\Gamma(m+1)}=\begin{cases}0&m=\dots,-2,-1;\\\frac1{m!},&m=0,1,\dots,\end{cases}
\eeq
we obtain
\begin{align}  {\mathbf F}(\sigma,\kappa,\omega;z)=
\sum_{n=\max\{0,-m\}}^\infty 
\frac{\Pi_{j=0}^{n-1}(\omega+(j+\12)\kappa'+\frac{j(j+1)}{2}\sigma'')}
  {(n+m)!n!}(-z)^n.\label{olver1}
  \end{align}
This easily implies the identity
\begin{align}
{\mathbf F}(\sigma,\kappa^\circ,\omega^\circ;z)=&\prod_{j=0}^{m-1}\Big(\omega-\kappa'(j+\tfrac12)+\frac{\sigma''}{2}j(j+1)\Big)(-z)^m 
                                                   {\mathbf F}(\sigma,\kappa,\omega;z), \quad m\in\nn_0.
\end{align}  
Therefore, for integer $m$  the functions  $\mathbf{F}(\sigma,\kappa,\omega;z)$ and $z^{-m}{\mathbf
  F}(\sigma,\kappa^\circ,\omega^\circ;z)$ are proportional to one
another and do
not form a basis of solutions.

\subsection{The power-exponential function}
\label{a.7}
In order to describe the theory of the  degenerate case in a unified
way, it will be convenient first  to unify the power  and  exponential
function in a single function. More precisely, let $a,\mu\in\cc$. 
Consider the function
\[
f(a,\mu;u):=\left\{\begin{array}{ll}
(1+\mu u)^{a\slash\mu},&\ \ \mu\neq0,\\[3mm]
\e^{au},&\ \ \mu=0.\end{array}\right.\]
It is analytic for $u\neq -\mu^{-1}$. From now on
we will write
$(1+\mu u)^{a\slash\mu}$ instead of $f(a,\mu;u)$.

Let us list some properties of this function.

\bet (1)
$(1+\mu u)^{a\slash\mu}$ is the unique solution of
\begin{align}
\big((1+\mu u)\p_u-a\big)f(u)=0,\\
f(0)=1.
\label{s1s}\end{align}
(2)
For $|u|<|\mu|^{-1}$ we have
\[(1+\mu u)^{a\slash\mu}
=\sum_{n=0}^\infty\frac{a(a-\mu)\cdots(a-(n-1)\mu)}{n!}u^n,
\]
(3)
\[(1+\mu u)^{a_1\slash\mu}(1+\mu u)^{a_2\slash\mu}
=(1+\mu u)^{(a_1+a_2)\slash\mu}.\]
\label{uni}\eet

\subsection{Generating functions}
\label{s.11}

Let us now fix $\sigma(z)=\frac{\sigma''}{2} z^2+z$, as in
\eqref{param}.
{\red Let $a,b,\mu,\nu\in\cc$. 
Consider two families of hypergeometric class operators
\begin{align}
  \Big(\frac{\sigma''}{2}
  z^2+z\Big)\partial_z^2+\Big(m+1+\Big(\frac{\sigma''}{2}(1+m)-\mu
  b-\nu a\Big)z\Big)\p_z-m\mu b-ab,\label{pasa1}\\
\Big(\frac{\sigma''}{2}
  z^2+z\Big)\partial_z^2+\Big(m+1+\Big(\frac{\sigma''}{2}(1+m)-\mu
  b-\nu a\Big)z\Big)\p_z-m\nu a-ab.\label{pasa2}
\end{align}
(Note that the parameters we introduced are redundant---three would
suffice instead of four. However they help us to describe the degenerate case in a
symmetric and unified way).

In terms of our standard parameters, the operators \eqref{pasa1} and
\eqref{pasa2} can be written as
\begin{align}
  {\red \cH}(\sigma,\kappa_m)+\omega_m,&\\
  {\red \cH}(\sigma,\kappa_m)+\tilde\omega_m,&
\end{align}
where
\begin{align}
  \kappa_m(z)&:=m+\Big(\frac{\sigma''}{2}(m-1)-\mu b-\nu a\Big)z,\\
  \omega_m&:=\frac{1}{2}(\mu b+\nu a)-ab-m \mu b-\frac{\sigma''}{4}(m-1),\\
  \tilde\omega_m&:=\frac{1}{2}(\mu b+\nu a)-ab-m \nu a-\frac{\sigma''}{4}(m-1).
\end{align}}
Note that
\beq\kappa_m^\circ=\kappa_{-m},\quad\omega_m^\circ=\tilde\omega_{-m},\quad\tilde\omega_m^\circ=\omega_{-m}.\eeq

By $[(\alpha,\beta)^+]$ we will denote a counterclockwise loop
around
 $\alpha,\beta\in\cc$.
We define
\begin{align}\label{psi1}
\Psi_m(z)&:=\frac{1}{2\pi\i}\int_{[(0,-z\nu))^+]}
(1+\mu u)^{-\frac{a}{\mu}}\Big(1+\frac{\nu z}{u}\Big)^{-\frac{b}{\nu}}u^{-m-1}\d u,\\\label{psi2}
\tilde\Psi_m(z)&:=\frac{1}{2\pi\i}\int_{[(0,-z\mu)^+]}
\Big(1+\frac{\mu z}{v}\Big)^{-\frac{a}{\mu}}(1+\nu v)^{-\frac{b}{\nu}}v^{-m-1}\d v.
\end{align}

\bet
 \begin{enumerate}\item
For
 $|z\nu|<|u|<|\mu^{-1}|$ we have
\beq
(1+\mu u)^{-\frac{a}{\mu}}\Big(1+\frac{\nu z}{u}\Big)^{-\frac{b}{\nu}}
=\sum_{m\in\zz}u^m\Psi_{m}(z).
\label{s3a}\eeq
 and for
 $|z\mu|<|u|<|\nu^{-1}|$ we have
\beq
\Big(1+\frac{\mu z}{v}\Big)^{-\frac{a}{\mu}}(1+\nu v)^{-\frac{b}{\nu}}
=\sum_{m\in\zz}v^m\tilde\Psi_{m}(z).
\label{s3ab}\eeq
\item
$
\Psi_{m}(z)=z^{-m}\tilde\Psi_{-m}(z).
$
\item
  For $m\geq0$,
\begin{align}
\Psi_m(0)&=
(-1)^m
\frac{a(a+\mu)\cdots(a+\mu(m-1))}{m!},\\
\tilde\Psi_m(0)&=(-1)^m\frac{b(b+\nu)\cdots(b+(m-1)\nu)}
{m!}.\end{align}
\item The functions $\Psi_{m}(z)$ and $\tilde\Psi_m(z)$ satisfy the
following hypergeometric class differential equations:
\begin{align}\left(
z(1-\mu\nu z)\p_z^2+
(m+1-(\mu\nu(1+m)+a\nu+b\mu)z)\p_z-(m\mu b+ab)\right)
\Psi_{m}(z)&=0,
             \label{s5}\\
  \left(
z(1-\mu\nu z)\p_z^2+
(m+1-(\mu\nu(1+m)+a\nu+b\mu)z)\p_z-(m\nu a+ab)\right)
\tilde\Psi_{m}(z)&=0.
\label{s5'}\end{align}
\item The functions $\Psi_{m}(z)$ and $\tilde\Psi_m(z)$ are
  proportional to {\red the Olver normalized unified} hypergeometric function:
\begin{align}
\Psi_m(z)&=
(-1)^m
a(a+\mu)\cdots(a+\mu(m-1)) \mathbf{F}(\sigma,\kappa_m,\omega_m;z);\\
  \tilde\Psi_m(z)&=(-1)^m b(b+\nu)\cdots(b+(m-1)\nu)
  \mathbf{F}(\sigma,\kappa_m,\tilde\omega_m;z);\end{align}
\end{enumerate}\eet

\proof 1. follows immediately from the definitions \eqref{psi1}, \eqref{psi2} and the Laurent
expansion. 

To show 2. we can rewrite
 (\ref{s3a}) and (\ref{s3ab}) as
\begin{align}
(1+\mu u)^{-\frac{a}{\mu}}(1+\nu v)^{-\frac{b}{\nu}}
&=\sum_{m\in\zz}u^m\Psi_{m}(uv),
\label{s6}\\
(1+\mu u)^{-\frac{a}{\mu}}(1+\nu v)^{-\frac{b}{\nu}}
&=\sum_{m\in\zz}v^m
\tilde\Psi_{m}(uv)=\sum_{m\in\zz}u^{-m}
(uv)^m\tilde\Psi_{m}(uv).
\label{s7}\end{align}

 {\red 3. follows by setting $v=0 $ in \eqref{s6} and $u=0$ in \eqref{s7}.}

Let us show 4. We have the identity
\begin{align}
0&=\left(\p_u\p_v-\mu\nu uv\p_u\p_v-\mu bu\p_u-a\nu v\p_v-ab\right)
(1+\mu u)^{-\frac{a}{\mu}}(1+\nu v)^{-\frac{b}{\nu}}.
\label{s8}\end{align}
Besides, 
\[\begin{array}{l}
\left(\p_u\p_v-\mu\nu uv\p_u\p_v-\mu bu\p_u-a\nu v\p_v-ab\right)
u^m\Psi_{m}(uv)\\[2mm]
=u^m\left(z(1-\mu\nu z)\p_z^2+
(m+1-(\mu\nu+(m\mu+a)\nu+b\mu)z)\p_z-(m\mu+a)b\right)
\Psi_{m}(z).\end{array}\]
Thus \eqref{s8} together with \eqref{s6} can be rewritten as
\beq
0=\sum_{m\in\zz}
u^m\left(z(1-\mu\nu z)\p_z^2+
(m+1-(\mu\nu+(m\mu+a)\nu+b\mu)z)\p_z-(m\mu+a)b\right)
\Psi_{m}(z),\eeq
which implies 
(\ref{s5}).

5. follows from 3. and 4.
\qed

In the remaining part of the section we describe separately  three
types of degenerate hypergeometric class functions.

\subsection{The ${}_2F_1$ function}

For $m\in\zz$ we have
\begin{align}{\bf F}
(a,b;1+m;z)&=\sum_{n=\max(0,-m)}\frac{(a)_n(b)_n}{n!(m+n)!}z^n,\\
(a-m)_m(b-m)_m{\bf F}(a,b;1+m;z)&=z^{-m}
{\bf F}(a-m,b-m;1-m;z).
\label{pwr}\end{align}

We have an integral representation
and a generating function:
\begin{align}
\frac{1}{2\pi \i}\int\limits_{[(0,z)^+]}(1-t)^{-a}\Big(1-\frac{z}{ t}\Big)^{-b}
t^{-m-1}\d t&=
(a)_m {\bf F}(a+m,b;1+m;z),\\
(1-t)^{-a}\Big(1-\frac{z}{ t}\Big)^{-b}
&=\sum_{m\in\zz}t^m(a)_m {\bf F} (a+m,b;1+m).
\end{align}

\subsection{The ${}_1F_1$ function}
If $m\in\zz$, we have
\begin{align}{\bf F}
(a;1+m;z)&=\sum_{n=\max(0,-m)}\frac{(a)_n}{n!(m+n)!}z^n,\\
(a-m)_m{\bf F}(a;1+m;z)&=z^{-m}
{\bf F}(a-m;1-m;z).
\label{pwrx}\end{align}

We have two
integral representations and the corresponding  generating functions:
\begin{eqnarray*}\frac{1}{2\pi\i}\int\limits
_{[(z,0)^+]}\e^{t}\Big(1-\frac{z}{t}\Big)^{-a}
t^{-m-1}\d t&=&
 {\bf F}(a,1+m;z),\\[4ex]
\frac{1}{2\pi\i}\int\limits_{[(0,1)^+]}
\e^{\frac{z}{ t}}(1-t)^{-a}t^{-m-1}\d t
&=&
 z^{-m} {\bf F} (a,1-m;z),\\
\e^{t}\Big(1-\frac{z}{ t}\Big)^{-a}&=&
\sum_{m\in\zz}t^m {\bf F} (a;m;z)),\\
\e^{\frac{z}{t}}(1-t)^{-a}&=&
\sum_{m\in\zz}t^m z^{-m}{\bf F}(a;1-m;z).\end{eqnarray*}

\subsection{The ${}_0F_1$ function}

If
$m\in\zz$, then
\begin{align}{\bf F}
(1+m;z)&=\sum_{n=\max(0,-m)}\frac{1}{n!(m+n)!}z^n,\\
{\bf F}(1+m;z)&=z^{-m}
{\bf F}(1-m;z).\end{align}

We have an integral representation, called the {\em Bessel formula},  and a generating function:
\begin{align}
\frac{1}{2\pi \i}\int\limits_{[0^+]}
  \e^{t+\frac{z}{t}}t^{-m-1}\d t&= {\bf F}  (1+m;z),
  \\
\e^{t}\e^{\frac{z}{t}}&
=\sum_{m\in\zz}t^m {\bf F} (1+m;z).
\end{align}

{\bf Acknowledgement.}
This  work 
 was supported by National Science Center (Poland) under the
    grant UMO-2019/35/B/ST1/01651.
I thank Rupert Frank for an inspiring remark.
I am also grateful to Adam Kali\'{n}ski for useful comments on the manuscript.
I have also greatly profited from many helpful remarks of the referee
of the previous version of the manuscript.

\appendix

{\red
\section{Function ${}_2F_0$}
\label{appendix}
\init

The function ${}_2F_0$ seems to be rarely discussed in the
literature. For convenience of the reader in the following theorem
we state and prove basic
facts about this function.

\bet For $w\in\cc\backslash[0,\infty[$, there exists the limit
\beq
F(a,b;-;w):=\lim_{c\to+\infty}F(a,b;c;cw).\eeq
It defines a function analytically depending on $a,b\in\cc$ and
$w\in\cc\backslash[0,\infty[$. 

We have the following asymptotic expansion:
\beq
F(a,b;-;w)\sim\sum_{n=0}^\infty
\frac{(a)_n(b)_n}{n!}w^n.\label{hehe1}\eeq
More precisely, for any $\epsilon>0$, $n$,
there exists $c_n$ such that
\beq
\Big|F(a,b;-;w)-\sum_{j=0}^n
\frac{(a)_j(b)_j}{j!}w^j\Big|\leq c_n|w|^{n+1},\quad
 |\arg w|\geq\epsilon,\quad |w|<1.\eeq
Moreover, for $\Re(a)>0$ we have an integral representation
\beq
F(a,b;-;w)=
\frac{1}{\Gamma(a)}\int_0^\infty\e^{-t}t^{a-1}(1-wt)^{-b}\d t.\label{f22}\eeq
\eet

\proof  Assume first that $\Re(a)>0$.
For $\Re(c-a)>0$ and $w\in\cc\backslash[c^{-1},+\infty[$, inserting
$t=cs^{-1}$ and $z=cw$ into \eqref{eqa1}, we obtain
\beq F(a,b;c;cw)=
\frac{\Gamma(c)c^{-a}}{\Gamma(a)\Gamma(c-a)}\int_0^c s^{a-1}(1-c^{-1}s
)^{c-a-1}(1-ws)^{-b}\d s.\label{hyhy}\eeq
Using $\lim_{c\to\infty}\frac{\Gamma(c)c^{-a}}{\Gamma(c-a)}=1$ and the
Lebesgue Dominated Convergence Theorem we see that
\eqref{hyhy} converges to the right hand side of \eqref{f22}. This
proves that $F(a,b;-;w)$ is well defined for $\Re(a)>0$.

Now let $a$ be arbitrary. We have
\beq \partial_z^nF(a,b;c;z)=\frac{(a)_n(b)_n}{(c)_n}
F(a+n,b+n;c+n;z).\eeq
Using Taylor's formula with a remainder 
\[f(z)=\sum_{j=0}^{n-1}\frac{f^{(j)}(0)
  z^j}{j!}+z^n\int_0^1\frac{f^{(n)}(sz)n(1-s)^{n-1}} {n!}\d s,\]
we obtain 
\begin{align}
F(a,b;c;z)
&=\sum_{j=0}^{n-1}
\frac{(b)_j(a)_j w^j}{j!}\\
&+
\frac{w^n(b)_n(a)_n}{(n-1)!}\int_0^1 (1-s)^{n-1} \d s
F(a+n,b+n;c;ws).
\end{align}
Now, choose $n$ such that  $\Re(a+n)>0$. Then we can apply what we
proved before to
show the convergence
\begin{align}
\lim_{c\to+\infty}F(a,b;c;cz)
&=\sum_{j=0}^{n-1}
\frac{(b)_j(a)_j w^j}{j!}\\
&+
\frac{w^n(b)_n(a)_n}{(n-1)!}\int_0^1 (1-s)^{n-1} \d s
F(a+n,b+n;-;ws).\label{psd1}
\end{align}
Clearly, \eqref{psd1} is $O(w^n)$, and $n$ can be made arbitrarily large.
\qed

If $a$ or $b$ is a negative integer, then the series \eqref{hehe1}
terminates and  functions ${}_2F_0(a,b;-;\cdot)$ essentially coincide with Bessel polynomials, see
Subsect. \ref{Bessel polynomials}. Otherwise, ${}_2F_0(a,b;-;\cdot)$
have a  logarithmic branch point at $0$.

\section{Riemann class}
\label{Class}
\init

We hope that we convinced our reader that the hypergeometric class
 studied in this paper is a natural family of equations. In this
appendix, following \cite{DIL},  we
 will describe a somewhat wider family, which  appears in the
 literature in many sources, e.g. \cite{SL,DIL,EEKS}.
According to the  terminology used in
\cite{DIL}, directly inspired by 
and partly borrowed from 
the monograph of
Slavyanov-Lay \cite{SL}, this family is called {\em Riemann class}.

Consider an equation given by the operator
  \begin{align}
&    \partial_z^2+b(z)\partial_z+c(z),
    \label{pde-}  \end{align}
  where $p(z)$, $q(z)$ are rational functions.
  Let $z_0\in\cc\cup\{\infty\}$ be a singular point of \eqref{pde-}.

  Recall that  $z_0\in\cc$ is called
{\em regular-singular} or {\em Fuchsian} if $b(z)=\frac{p(z)}{z-z_0}$ and
$c(z)=\frac{q(z)}{(z-z_0)^2}$ with $p,q$
 regular at $z_0$. 
The equation
\beq 
\lambda(\lambda-1)+p(z_0)\lambda+q(z_0)=0 
\eeq 
is called the {\em indicial equation of $z_0$} and its roots are called 
{\em indices} of $z_0$.

We say that $\infty $ is 
{\em regular-singular} or {\em Fuchsian} if $b(z)=\frac{p(z)}{z}$ and 
$c(z)=\frac{q(z)}{z^2}$ with $p,q$ regular 
at $\infty$.
\beq 
\lambda(\lambda+1)-p(\infty)\lambda+q(\infty)=0 
\eeq 
is called the {\em indicial equation of $\infty$} and its roots are called 
{\em indices} of $\infty$.

It is easy to see that every  equation having no more than $3$  singular points in $\cc\cup\{\infty\}$,  all of them Fuchsian
and at most $2$  finite, is given by an operator of the form
\begin{equation}
  \partial_z^2+\Big(\frac{a_1}{z-z_1}+\frac{a_2}{z-z_2}\Big)\partial_z
  +\frac{b_1}{z-z_1}  +\frac{b_2}{z-z_2}
    +\frac{c_1}{(z-z_1)^2}  +\frac{c_2}{(z-z_2)^2}
  \label{req1-}
\quad  \text{ with  }\quad
b_1+b_2=0,\end{equation}
where $z_1,z_{2}$ are distinct points in $\cc$.
Following \cite{DIL}, the family of equations (\ref{req1-}) is called
{\em the Riemann type}. (Another name, introduced in \cite{SL}, is $M_2$-type).

Each finite singularity has at least one index equal $0$ if and only if
$c_1=c_{2}=0.$
Such equations are given by operators 
\begin{equation}
  \partial_z^2+\Big(\frac{a_1}{z-z_1}+\frac{a_2}{z-z_2}\Big)\partial_z
  +\frac{b_1}{z-z_1}  +\frac{b_2}{z-z_2}
\label{req2}
\quad  \text{ with  }\quad
  \sum_{j=1}^{n}b_j=0.\end{equation}
Following \cite{DIL}, the family of equations given by (\ref{req2}) will be called
{\em  the grounded Riemann  type}.
It is easy to see that by gauging with power functions we can always transform a
  Riemann
  type equation into a grounded Riemann  type equation.

We say that a differential equation belongs to the {\em Riemann
  class} (or the $M_2$-class) if it is given by
\begin{equation}
  \partial_z^2+\frac{\tau(z)}{\sigma(z)}\partial_z+\frac{\xi(z)}{\sigma(z)^2},
  \label{req3}\end{equation}
where $\sigma,\tau,\xi$ are polynomials satisfying
\begin{align}\label{req3/}
  \sigma\neq0,\quad \deg\sigma\leq 2,\quad\deg\tau\leq 1,\quad
\deg\xi\leq2.
\end{align}
Thus the Riemann class comprises the Riemann type together with all
its confluent cases.

We say that a differential equation belongs to the {\em grounded Riemann class} if it is given by
\begin{equation}
  \partial_z^2+\frac{\tau(z)}{\sigma(z)}\partial_z+\frac{\eta(z)}{\sigma(z)},
  \label{req3-}\end{equation}
where $\sigma,\tau,\eta$ are polynomials satisfying
\begin{align}
    \label{req3-/}
  \sigma\neq0,\quad\deg\sigma\leq 2,\quad\deg\tau\leq 1,\quad
\deg\eta=0.
\end{align}
Thus the grounded Riemann class besides the grounded Riemann includes all 
its confluent cases.

The following  proposition is proven in \cite{DIL}:
\begin{proposition}\begin{enumerate}\item
The Riemann type is contained in the Riemann class. An equation of the
Riemann class  is of the Riemann type iff $\sigma$ possesses $2$
distinct roots. 
\item
The grounded Riemann type is the intersection of the grounded Riemann
class and Riemann type.
\end{enumerate}  \end{proposition}

In this paper instead of the {\em grounded Riemann class} we use the
name
{\em hypergeometric class}.
We use this name mostly for brevity, besides the word ``grounded'' may
sound bizarre to some readers.
Furthermore, we prefer to multiply these
equations by $\sigma(z)$, so that we consider operators of the form
\begin{equation}
 \sigma(z) \partial_z^2+\tau(z)\partial_z+\eta
  \label{req3.}\end{equation}
with $\sigma,\tau,\nu$ satisfying \eqref{req3-/}.

One can ask whether the unified theory
presented in this paper can be extended to operators of the form
\begin{equation}
 \sigma(z) \partial_z^2+\tau(z)\partial_z+\frac{\xi(z)}{\sigma(z)}, 
  \label{req3}\end{equation}
where $\sigma,\tau,\xi$  satisfy
\eqref{req3/}. In other words whether the unified theory presented in
this paper can be extended to the full Riemann class.
We do not know an answer to this question. However, let us remark that
we do not gain much by considering the full Riemann class instead of
the grounded Riemann class.

By a division by a constant,
transformation $z\mapsto az+b$ and gauging with a power and
exponential  all
operators from the Riemann type
 can be always reduced
 to an ${}_2F_1$ operator.
 More generally, all operators from the Riemann class can be reduced
 to one of normal forms listed in \cite{SL} or  in Subsect. 3.3 of
 \cite{DIL}:
\medskip
 
\begin{tabular}{lll}
  \hline\\
the  ${}_2F_1$ operator&
                                                $z(1-z)\p_z^2+\big(c-(a+b+1)z\big)\p_z-ab$;\\[2ex]
  the  ${}_2F_0$ operator& 
                                                 $z^2\p_z^2+\big(-1+(a+b+1)z\big)\p_z+ab$;\\[2ex]
  the   ${}_1F_1$ operator& 
                                                  $z\p_z^2+(c-z)\p_z-a$;\\[2ex]
  the   ${}_0F_1$ operator& 
                                                        $z\p_z^2+c\p_z-1$;\\[2ex]
  the Hermite operator&
                                  $\p_z^2-2z\p_z-2a$;\\[2ex]
  the Airy operator &
                                     $\p_z^2+z$;\\[2ex]
  the Euler  operator I& 
                                                $z^2\p_z^2+cz\p_z$;\\[2ex]
  the Euler  operator II& 
                                                $z\p_z^2+c\p_z$;\\[2ex]
  the 1d Helmholtz operator &
                                       $\p_z^2+1$;\\[2ex]
  the 1d Laplace operator & 
                                                 $\p_z^2$.\\
 \hline
\end{tabular}
\medskip

In this list the first five 
 are the normal forms listed in the introduction and
described in 
Section \ref{Applications case by case}.
The last four are trivial operators that can be solved in by
elementary methods. They are also  special cases of the
grounded Riemann class.

The Airy equation is the only type within  the Riemann class
which does not belong to the grounded Riemann class and therefore is 
left out from the unified theory presented in this paper.

}

\end{document}